# Global Stabilization of Compressible Flow Between Two Moving Pistons

**Iasson Karafyllis[*] and Miroslav Krstic[**]**

[*]Dept. of Mathematics, National Technical University of Athens,
Zografou Campus, 15780, Athens, Greece, email: iasonkar@central.ntua.gr

[**]Dept. of Mechanical and Aerospace Eng., University of California, San Diego, La Jolla, CA 92093-0411, U.S.A., email: krstic@ucsd.edu

**Abstract**

This paper studies the global feedback stabilization problem of a system with two pistons and the area between them containing a viscous compressible fluid (gas) modeled by the Navier-Stokes equations. The control input is the force applied on the left piston (boundary input) and the overall system consists of two nonlinear Partial Differential Equations and four nonlinear Ordinary Differential Equations. Global feedback stabilizers are designed for the overall system by means of the Control Lyapunov Functional methodology. The closed-loop system exhibits global asymptotic stability with an exponential convergence rate. The proposed stabilizing boundary feedback laws do not require measurement of the density and velocity profiles inside the area between the pistons and simply require measurements of the gas density and velocity at the position of the actuated piston.

## 1. Introduction

In recent years, the interest has risen in modeling, stability analysis and control of systems containing compressible fluids moving in one dimension (see for instance [1,2,3,7,8,9,10,14]). Compressible fluids are modeled by means of the Navier-Stokes equations and their mathematical analysis requires special tools to deal with the nonlinearity (see for instance [11,17,20,25,30]). Among the problems with compressible fluids, the problems that describe the 1-D movement of pistons are of considerable industrial interest. The Partial Differential Equations (PDEs) that are used in the description of the movement of a piston may be divided into two main categories:

1) Models that take into account the energy balance (see [12,29,30]). In these models, each compartment containing a compressible gas is modeled by three PDEs (the mass conservation equation, the momentum conservation equation and the energy conservation equation) and the three states are the density, the velocity and the temperature of the fluid.

2) Models that do not take into account the energy balance but exploit a relationship between pressure and density based on the assumption of isentropic or barotropic flow (see [4,13,23,26,27,28,30]). In these models, each compartment containing a compressible gas is modeled by two PDEs (the mass conservation equation and the momentum conservation equation) and the two states are the density and the velocity of the fluid.

There are many other categorizations that may be used for the piston problems: whether the gas is viscous or inviscid (giving rise to systems of hyperbolic PDEs in the latter case and systems of hyperbolic-parabolic PDEs in the former case) or whether the PDEs are given in Lagrange coordinates or Euler coordinates.

One of the great difficulties that arise in the analysis of piston problems is the fact that the movement of the piston introduces a moving boundary. The 1-D movement of the piston itself is modeled by a pair of Ordinary Differential Equations (ODEs) that express Newton's second law. This feature becomes a major challenge when control problems for a piston are studied. Recently, control problems for systems with moving boundaries have attracted the attention of researchers; see [18] and references therein for control problems that deal with phase changes (Stefan problems). Another challenge that appears in piston problems is the positivity constraint for the density (negative values for the density are physically irrelevant). In addition to the above, the nonlinear character of the Navier-Stokes equations is another major difficulty. The number of difficulties and challenges that appear in a piston control problem explains the scarcity of the available results in the literature. There are very few results in the literature for the control of pistons:

a) in [22] the authors study the null-controllability of a simplified piston model with one input, and

b) in [26] the authors formulate an optimal control problem for a piston with multiple boundary inputs.

The only feedback control design for a piston control problem was performed in [4], where the authors used the linearized model of an inviscid gas with two boundary inputs (the gas velocities on the boundaries).

This paper focuses on the global feedback stabilization problem of a system with two pistons. The area between the pistons contains a viscous compressible fluid (gas), which is modeled by the Navier-Stokes equations. There is only one boundary control input: an external force acting on the left piston. We do not apply any linearization and we consider the nonlinear model. The nonlinear model implies an algebraic constraint (the conservation of mass) that must also be satisfied. In order to make the model more realistic, we take into account the (nonlinear) dependence of the viscosity on the fluid temperature and density and we assume a very general relationship between the pressure and the density of the gas. Global feedback stabilizers are designed for the overall system that contains two PDEs and four ODEs. Global asymptotic stability with an exponential convergence rate is achieved for the closed-loop system. The stabilizing feedback laws do not require measurement of the density and velocity spatial profiles between the pistons and simply require measurements of the gas density and velocity at the position of the actuated piston. Moreover, if certain additional assumptions about the physical properties of the gas are fulfilled, then it is not even necessary to measure the boundary gas density and it is not necessary to know the relation between density and pressure and the relation between viscosity and density. To our knowledge, this is the first paper in the literature that achieves global stabilization of the two-piston system. In addition, we show that all our assumptions are valid for an ideal gas under constant entropy but are also valid for a wide class of non-ideal gases.

The boundary feedback design methodology that is followed in the present work is a Control Lyapunov Functional (CLF) methodology that was first used for global stabilization of nonlinear parabolic PDEs in [19] and was subsequently studied in [15,16] (but see also [6] for the presentation of the CLF methodology in finite-dimensional and infinite-dimensional systems). The CLF plays also the role of a barrier function and guarantees a positive lower bound for the density. The CLF is constructed by combining the mechanical energy of the system and the use of a specific transformation that has been used in the literature of compressible fluids (see [23,27,30]). It is important to note here that the nonlinearity of the control problem hinders the application of standard boundary feedback design methodologies like backstepping (see [31,32,34]; although backstepping has been applied successfully for problems with moving boundary-see [18]).

The structure of the paper is as follows. Section 2 is devoted to the presentation of the control problem. Section 3 contains the statements of the assumptions and the main results of the paper. Moreover, in Section 3 we also provide insights for the proofs of the main results: we explain the construction of the CLF, the selection of the state space and we provide some auxiliary results that are used in the proofs. The proofs of all results are provided in Section 4. Finally, Section 5 gives the concluding remarks of the present work.

**Notation.** Throughout this paper, we adopt the following notation.

* $\Re_+ = [0, +\infty)$ denotes the set of non-negative real numbers.

* Let $S \subseteq \Re^n$ be an open set and let $A \subseteq \Re^n$ be a set that satisfies $S \subseteq A \subseteq cl(S)$. By $C^0(A; \Omega)$, we denote the class of continuous functions on $A$, which take values in $\Omega \subseteq \Re^m$. By $C^k(A; \Omega)$, where $k \geq 1$ is an integer, we denote the class of functions on $A \subseteq \Re^n$, which takes values in $\Omega \subseteq \Re^m$ and has continuous derivatives of order $k$. In other words, the functions of class $C^k(A; \Omega)$ are the functions which have continuous derivatives of order $k$ in $S = \text{int}(A)$ that can be continued continuously to all points in $\partial S \cap A$. When $\Omega = \Re$ then we write $C^0(A)$ or $C^k(A)$. When $I \subseteq \Re$ is an interval and $\eta \in C^1(I)$ is a function of a single variable, $\eta'(\rho)$ denotes the derivative with respect to $\rho \in I$.

* Let $I \subseteq \Re$ be an interval and let $a, b : I \to \Re$ with $a(t) < b(t)$ for all $t \in I$, $u : \left( \bigcup_{t \in I} (\{t\} \times [a(t), b(t)]) \right) \to \Re$ be given. We use the notation $u[t]$ to denote the profile at certain $t \in I$, i.e., $(u[t])(x) = u(t, x)$ for all $x \in [a(t), b(t)]$. When $u(t, x)$ is (twice) differentiable with respect to $x \in [a(t), b(t)]$, we use the notation $u_x(t, x)$ ($u_{xx}(t, x)$) for the (second) derivative of $u$ with respect to $x \in [a(t), b(t)]$, i.e., $u_x(t, x) = \frac{\partial u}{\partial x}(t, x)$ ($u_{xx}(t, x) = \frac{\partial^2 u}{\partial x^2}(t, x)$). When $u(t, x)$ is differentiable with respect to $t$, we use the notation $u_t(t, x)$ for the derivative of $u$ with respect to $t$, i.e., $u_t(t, x) = \frac{\partial u}{\partial t}(t, x)$.

* Given a set $U \subseteq \Re^n$, $\chi_U$ denotes the characteristic function of $U$, i.e. the function defined by $\chi_U(x) := 1$ for all $x \in U$ and $\chi_U(x) := 0$ for all $x \notin U$.

* Let $a < b$ be given constants. For $p \in [1, +\infty)$, $L^p(a, b)$ is the set of equivalence classes of Lebesgue measurable functions $u : (a, b) \to \Re$ with $\|u\|_p := \left( \int_a^b |u(x)|^p dx \right)^{1/p} < +\infty$. $L^\infty(a, b)$ is the set of equivalence classes of Lebesgue measurable functions $u : (a, b) \to \Re$ with $\|u\|_\infty := \text{ess} \sup_{x \in (a,b)} (|u(x)|) < +\infty$. For an integer $k \geq 1$, $H^k(a, b)$ denotes the Sobolev space of functions in $L^2(a, b)$ with all its weak derivatives up to order $k \geq 1$ in $L^2(a, b)$.

## 2. Description of the Problem

We consider a one dimensional model for the motion of two pistons. The area between the pistons contains a viscous compressible gas. The right side of the right piston is subject to a constant external pressure, while the left side of the left piston is subject to a pressure-actuated force that can be manipulated. No gas is allowed to penetrate through the two pistons. The gas is modeled by the 1-D compressible Navier–Stokes equations for isentropic motion, whereas the pistons obey Newton's second law. Since the positions of the pistons (and consequently, the domain occupied by the gas) are unknowns of the problem, we have a free boundary problem. Figure 1 shows a picture of the two-piston system.

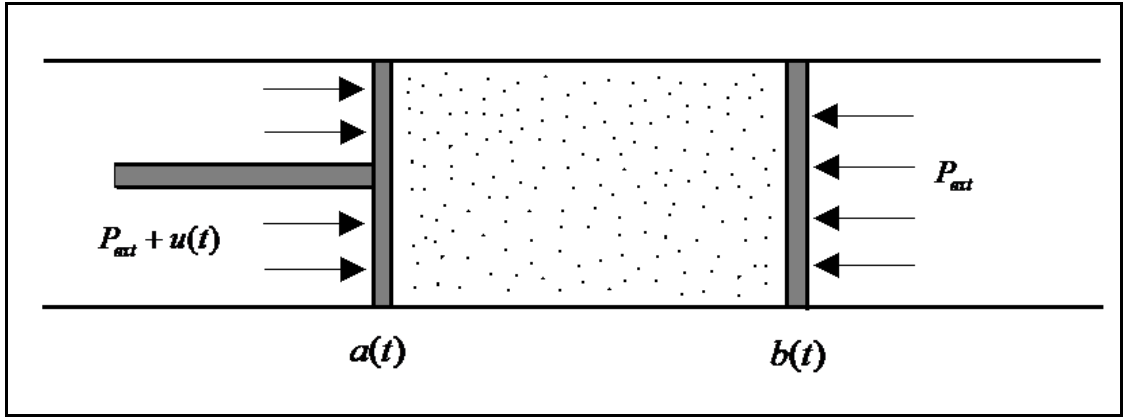

**Figure 1:** The two-piston system.

Let the position of the left piston at time $t \geq 0$ be $a(t)$ and let the position of the right piston at time $t \geq 0$ be $b(t) > a(t)$. The equations describing the motion of the isentropic (or barotropic; see [33], page 107) viscous gas between the two pistons are

$$\frac{\partial \rho}{\partial t}(t,x) + \frac{\partial}{\partial x}\big(\rho(t,x)v(t,x)\big) = 0, \text{ for } t > 0, \ x \in [a(t), b(t)] \tag{2.1}$$

$$\rho(t,x)\frac{\partial v}{\partial t}(t,x) + \rho(t,x)v(t,x)\frac{\partial v}{\partial x}(t,x) = \frac{\partial}{\partial x}\left(\mu(\rho(t,x))\frac{\partial v}{\partial x}(t,x) - P(\rho(t,x))\right)$$
$$\text{for } t > 0, \ x \in \big(a(t), b(t)\big) \tag{2.2}$$

where $\rho(t,x) > 0$, $v(t,x) \in \Re$ are the density and velocity of the gas, respectively, at time $t \geq 0$ and position $x \in [a(t), b(t)]$, while $P(\rho(t,x)), \mu(\rho(t,x)) > 0$ are the pressure and viscosity of the gas, respectively, at time $t \geq 0$ and position $x \in [a(t), b(t)]$. Both functions $P, \mu : (0, +\infty) \to (0, +\infty)$ are assumed to be $C^2$ functions of density and $P$ is also assumed to be increasing with $P\big((0, +\infty)\big) = (0, +\infty)$. Typically, the function $P(\rho) = c\rho^\gamma$, where $c, \gamma > 0$ are constants with $\gamma \in (1, 2)$ is the most common function used in the literature but here we will not specify the functions $P, \mu : (0, +\infty) \to (0, +\infty)$.

The velocities of the pistons coincide with the gas velocities on both sides, i.e., we have:

$$v(t, a(t)) = \dot{a}(t), \ v(t, b(t)) = \dot{b}(t), \text{ for } t \geq 0 \tag{2.3}$$

Moreover, applying Newton's second law for the two pistons we get

$$\ddot{b}(t) = P(\rho(t, b(t))) - \mu(\rho(t, b(t)))v_x(t, b(t)) - P_{ext}, \text{ for } t > 0 \tag{2.4}$$

where $P_{ext} > 0$ is the constant external pressure applied to the right side of the right piston, and

$$\ddot{a}(t) = u(t) + P_{ext} - P(\rho(t, a(t))) + \mu(\rho(t, a(t)))v_x(t, a(t)), \text{ for } t > 0 \quad (2.5)$$

where $P_{ext} + u(t)$ is the manipulated pressure applied to the left side of the left piston. The variable $u(t) \in \Re$ is the control input and is the deviation of the pressure applied to the left side of the left piston from the equilibrium pressure $P_{ext} > 0$.

We consider classical solutions for the PDE-ODE system (2.1)-(2.5), i.e., we consider functions $a, b \in C^1(\Re_+) \cap C^2((0, +\infty))$ satisfying $b(t) > a(t)$ for all $t \geq 0$, $\rho \in C^1(\bar{\Omega}; (0, +\infty)) \cap C^2(\Omega)$, $v \in C^0(\bar{\Omega}) \cap C^1(\tilde{\Omega})$ with $v[t] \in C^2((a(t), b(t)))$ for each $t > 0$ and $\Omega = \bigcup_{t>0}(\{t\} \times (a(t), b(t)))$, $\tilde{\Omega} = \bigcup_{t>0}(\{t\} \times [a(t), b(t)])$ that satisfy equations (2.1)-(2.5) for a given input $u \in C^0(\Re_+)$.

Using (2.1) and (2.3), we can prove that for every solution of (2.1)-(2.5) it holds that $\dot{m}(t) = 0$ for all $t > 0$, where $m(t)$ is the total mass of the gas

$$m(t) = \int_{a(t)}^{b(t)} \rho(t, x) dx$$

Therefore, the total mass of the gas is constant. Without loss of generality, we assume that every solution of (2.1)-(2.5) satisfies the equation

$$\int_{a(t)}^{b(t)} \rho(t, x) dx \equiv 1 \quad (2.6)$$

The open-loop system (2.1)-(2.6), i.e., system (2.1)-(2.6) with $u(t) \equiv 0$, allows a continuum of equilibrium points, namely the points

$$\rho(x) \equiv \rho^*, \ v(x) \equiv 0, \text{ for } x \in [a, b] \quad (2.7)$$

$$b = a + 1/\rho^* \quad (2.8)$$

where $\rho^* > 0$ is the unique solution of the equation

$$P(\rho^*) = P_{ext} \quad (2.9)$$

Notice that the equilibrium position of the left piston $a \in \Re$ is not determined by equations (2.7), (2.8), (2.9) and is a free parameter. Therefore, there are infinite equilibrium points for the open-loop system (2.1)-(2.6).

Our objective is to design a boundary feedback law of the form

$$u(t) = F(\rho(t, a(t)), \dot{a}(t)), \text{ for } t > 0 \quad (2.10)$$

which achieves global stabilization of the set of equilibrium points (without having the ability to prescribe a particular equilibrium point). The term "global" means that *every solution* of the closed-loop system (2.1)-(2.6), (2.10) satisfies a specific stability estimate that shows convergence to the desired set. The equilibrium values $\rho = \rho^*$ and $b - a = 1/\rho^*$ to which the closed-loop solution converges, will depend on the initial condition.

It should be emphasized that the feedback law (2.10) does not require measurement of the whole state $(a(t), b(t), \rho[t], v[t])$ (notice that $\dot{a}(t)$ and $\dot{b}(t)$ are included in the state vector because of the boundary condition (2.3)) but only measurements at the position of the actuated piston. This feature is very important for the implementation of the feedback law (2.10).

Finally, it should be noticed that the desired set of equilibrium points is not asymptotically stable for the open-loop system. The non-convergence of solutions of the open-loop system that start from points which are arbitrarily close to the desired equilibrium point can be shown by the preservation of the momentum for the open-loop system, i.e., from the fact that for $u(t) \equiv 0$ every solution of (2.1)-(2.6) satisfies

$$\frac{d}{dt}\left( \dot{a}(t) + \dot{b}(t) + \int_{a(t)}^{b(t)} \rho(t,x) v(t,x) dx \right) \equiv 0$$

Consequently, the total momentum $R(t) = \dot{a}(t) + \dot{b}(t) + \int_{a(t)}^{b(t)} \rho(t,x) v(t,x) dx$ is constant for every solution of the open-loop system. Hence only solutions of the open-loop system that satisfy $\dot{a}(t) + \dot{b}(t) + \int_{a(t)}^{b(t)} \rho(t,x) v(t,x) dx \equiv 0$ can converge to one of the equilibrium points described by (2.7), (2.8), (2.9).

## 3. Construction of the Boundary Feedback Law

3.1. Physical properties of the gas

For the study of the two-piston problem we need a technical assumption involving the physical properties of the gas.

**(H)** *The functions $P, \mu : (0, +\infty) \to (0, +\infty)$ are positive $C^2$ functions with $P'(\rho) > 0$ for all $\rho > 0$ and $P\big((0, +\infty)\big) = (0, +\infty)$. Moreover, for every $\rho^* > 0$ the following equations hold:*

$$\lim_{\rho \to +\infty} \big(G(\rho)\big) = +\infty \qquad (3.1)$$

$$\lim_{\rho \to 0^+} \big(G(\rho)\big) = -\infty \qquad (3.2)$$

$$\lim_{\rho \to +\infty} \big(k(\rho)\big) = +\infty \qquad (3.3)$$

*where $G, k : (0, +\infty) \to \Re$ are the functions*

$$G(\rho) := \int_{\rho^*}^{\rho} \frac{\mu(l)}{l^{3/2}} \sqrt{Q(l)} dl \qquad (3.4)$$

$$k(\rho) := \int_{\rho^*}^{\rho} \frac{\mu(\tau)}{\tau} d\tau \qquad (3.5)$$

*and $Q : (0, +\infty) \to \Re_+$ is the function*

$$Q(\rho) := \rho \int_{\rho^*}^{\rho} \frac{P(\tau)}{\tau^2} d\tau - \frac{P(\rho^*)}{\rho^*} \rho + P(\rho^*) \qquad (3.6)$$

**Remarks: (i)** Since $Q''(\rho) = \rho^{-1}P'(\rho) > 0$ for all $\rho > 0$, $Q(\rho^*) = Q'(\rho^*) = 0$, it follows that $Q(\rho) > 0$ for all $\rho > 0$, $\rho \neq \rho^*$.

**(ii)** It should be noticed that Assumption (H) is valid for an ideal gas under constant entropy. Indeed, for an ideal gas under constant entropy we have $P(\rho) = c\rho^\gamma$, $P = f\rho T$ and $\mu = \beta\sqrt{T}$, where $T > 0$ is the temperature and $c, f, \beta, \gamma > 0$ are constants with $\gamma \in (1,2)$ (see [24]). Therefore, we have $\mu(\rho) = A\rho^{\frac{\gamma-1}{2}}$ with $A = \beta\sqrt{f^{-1}c}$ and $Q(\rho) := \frac{c}{\gamma-1}\left(\rho^\gamma - \gamma(\rho^*)^{\gamma-1}\rho + (\gamma-1)(\rho^*)^\gamma\right)$, $k(\rho) := \frac{2A}{\gamma-1}\left(\rho^{\frac{\gamma-1}{2}} - (\rho^*)^{\frac{\gamma-1}{2}}\right)$, $G(\rho) := A\int_{\rho^*}^{\rho} l^{\frac{\gamma-4}{2}}\sqrt{Q(l)}dl$. Notice that the fact that $Q(\rho) \geq \frac{c}{2(\gamma-1)}\rho^\gamma$ for $\rho \geq \frac{\ln(2\gamma)}{\gamma-1}\rho^*$ and the fact that $Q(\rho) \geq \frac{c}{2}(\rho^*)^\gamma$ for $0 < \rho < \frac{\gamma-1}{2\gamma}\rho^*$ imply that

$$G(\rho) \geq A\sqrt{\frac{c}{2(\gamma-1)}} \frac{\rho^{\gamma-1} - \left(\frac{\ln(2\gamma)}{\gamma-1}\rho^*\right)^{\gamma-1}}{\gamma-1} \quad \text{for} \quad \rho \geq \frac{\ln(2\gamma)}{\gamma-1}\rho^* \quad \text{and}$$

$$G(\rho) \leq -\frac{2A}{2-\gamma}\sqrt{\frac{c}{2}}(\rho^*)^{\frac{\gamma}{2}}\left(\rho^{\frac{\gamma-2}{2}} - \left(\frac{\gamma-1}{2\gamma}\rho^*\right)^{\frac{\gamma-2}{2}}\right) \quad \text{for} \quad 0 < \rho < \frac{\gamma-1}{2\gamma}\rho^*. \text{ Consequently, the fact that}$$

$\gamma \in (1,2)$ guarantees that Assumption (H) is valid for an ideal gas under constant entropy.

**(iii)** Working as above it is possible to show that equations (3.1), (3.2), (3.3) are valid for any gas that satisfies the inequalities $P(\rho) \geq c\rho^\gamma$ and $\mu(\rho) \geq A\rho^\eta$, where $c, A, \gamma > 0$, $\eta \in \Re$ are constants with $\gamma \in (1,2)$ and $\eta \in [0, 1/2]$. Therefore, equations (3.1), (3.2), (3.3) are valid for a wide class of gases under constant entropy.

3.2. The Control Lyapunov Functional (CLF)

Let $r > 0$ be a constant (to be selected). Given $a,b \in \Re$ with $a < b$, $\rho \in C^0([a,b];(0,+\infty)) \cap H^1(a,b)$, $v \in C^0([a,b])$ with $\int_a^b \rho(x)dx = 1$, we define the following functionals:

$$V(a,b,\rho,v) := W(a,b,\rho,v) + rE(a,b,\rho,v) \tag{3.7}$$

$$E(a,b,\rho,v) := \frac{v^2(a)}{2} + \frac{v^2(b)}{2} + \frac{1}{2}\int_a^b \rho(x)v^2(x)dx + U(a,b,\rho) \tag{3.8}$$

$$W(a,b,\rho,v) := \frac{1}{2}\int_a^b \frac{1}{\rho(x)}\left(\rho(x)v(x) + \frac{\partial}{\partial x}(k(\rho(x)))\right)^2 dx \\ + \frac{1}{2}(v(b) - k(\rho(b)))^2 + \frac{1}{2}(v(a) + k(\rho(a)))^2 + U(a,b,\rho) \tag{3.9}$$

$$U(a,b,\rho) := \int_a^b Q(\rho(x))dx \tag{3.10}$$

We notice that:

- the functional $U$ is the potential energy of the two-piston system,
- the functional $E$ is the total mechanical energy of the two-piston system. Indeed, notice that $E$ is the sum of the potential energy ($U$) and the kinetic energy of the two-piston system: $v^2(a)/2$ is the kinetic energy of the left piston (recall (2.3)), $v^2(b)/2$ is the kinetic energy of the right piston (recall (2.3)) and $\frac{1}{2}\int_a^b \rho(x)v^2(x)dx$ is the kinetic energy of the gas,
- the functional $W$ is a kind of mechanical energy and has been constructed based on a specific transformation that has been used extensively in the literature of isentropic, compressible fluid flow (see [23,27,30]). More specifically, the transformation $w = \rho v + (k(\rho))_x$ has the effect of making the viscosity term disappear from the momentum equation, which then becomes $w_t + (vw)_x = -(P(\rho))_x$.

We use the functional $V(a,b,\rho,v)$ defined by (3.7) as a CLF for the two-piston system. However, the functional $V(a,b,\rho,v)$ plays also of the role of a barrier function in the sense that $V(a,b,\rho,v) \to +\infty$ as $\min_{a \leq x \leq b}(\rho(x)) \to 0$. This is guaranteed by the following lemma.

**Lemma 1:** *Suppose that Assumption (H) is valid. For every $S \geq 0$ there exist constants $\rho_{\max} \geq \rho_{\min} > 0$ (depending only on $S \geq 0$, $r > 0$ and $\rho^* > 0$) satisfying the following property:*

*"for every $a,b \in \Re$ with $a < b$, $\rho \in C^0([a,b];(0,+\infty)) \cap H^1(a,b)$, $v \in C^0([a,b])$ with $\int_a^b \rho(x)dx = 1$ and*

$$V(a,b,\rho,v) \leq S \tag{3.11}$$

*where $V$ is defined by (3.7), the following inequality holds:*

$$\rho_{\min} \leq \rho(x) \leq \rho_{\max}, \text{ for all } x \in [a,b]" \tag{3.12}$$

3.3. The state space

In order to be able to give the main results of the present work, we need to explain the structure of the state space. Given $a,b \in \Re$, $a < b$, $\rho \in C^0([a,b];(0,+\infty)) \cap H^1(a,b)$, $v \in C^0([a,b])$ with $\int_a^b \rho(x)dx = 1$, consider the operator $(a,b,\rho,v) \to \Phi(a,b,\rho,v) = (\tilde{\rho},\tilde{v})$ with

$$\tilde{\rho}(\theta) = \rho(a+(b-a)\theta), \tilde{v}(\theta) = v(a+(b-a)\theta), \text{ for } \theta \in [0,1] \tag{3.13}$$

Notice that $\Phi(a,b,\rho^*\chi_{[a,b]}, 0\chi_{[a,b]}) = (\rho^*\chi_{[0,1]}, 0\chi_{[0,1]})$ for all $a \in \Re$ and $b = a + (\rho^*)^{-1}$. Clearly, $\Phi(a,b,\rho,v) = (\tilde{\rho},\tilde{v})$ as defined by (3.13) takes values in the convex cone $\tilde{X} = (C^0([0,1];(0,+\infty)) \cap H^1(0,1)) \times C^0([0,1]) \subset X = H^1(0,1) \times C^0([0,1])$. We consider the linear space $X = H^1(0,1) \times C^0([0,1])$ endowed with the following norm for $(\rho,v) \in X$:

$$\|(\rho,v)\|_X = \sqrt{v^2(0)+v^2(1)} + \left(\int_0^1 v^2(\theta)d\theta\right)^{1/2} + \max_{0 \leq \theta \leq 1}(|\rho(\theta)|) + \left(\int_0^1 (\rho'(\theta))^2 d\theta\right)^{1/2} \tag{3.14}$$

Consider next a classical solution $(a(t),b(t),\rho[t],v[t])$ for $t\geq 0$ of (2.1)-(2.6), (2.10) with $a(t),b(t)\in\Re$, $a(t)<b(t)$, $\rho[t]\in C^0\left([a(t),b(t)];(0,+\infty)\right)\cap H^1(a(t),b(t))$, $v[t]\in C^0\left([a(t),b(t)]\right)$ with $\int_{a(t)}^{b(t)}\rho(t,x)dx=1$. We define:

$$(\tilde{\rho}[t],\tilde{v}[t])=\Phi(a(t),b(t),\rho[t],v[t])\in\tilde{X}, \text{ for } t\geq 0 \quad (3.15)$$

The mapping defined by (3.15) is "1-1" for a fixed $a(0)\in\Re$ in the sense that given $(\tilde{\rho}[t],\tilde{v}[t])\in\tilde{X}$ for $t\geq 0$ we are in a position to reconstruct the solution $(a(t),b(t),\rho[t],v[t])$ for $t\geq 0$ by means of the formulas:

$$a(t)=a(0)+\int_0^t \tilde{v}(l,0)dl$$

$$b(t)=a(t)+\left(\int_0^1 \tilde{\rho}(t,\theta)d\theta\right)^{-1} \quad (3.16)$$

$$\rho(t,x)=\tilde{\rho}\left(t,\frac{x-a(t)}{b(t)-a(t)}\right), v(t,x)=\tilde{v}\left(\frac{x-a(t)}{b(t)-a(t)}\right), \text{ for } x\in[a(t),b(t)] \quad (3.17)$$

Indeed, since (2.3) holds, we obtain from (3.13), (3.15) that $\dot{a}(t)=\tilde{v}(t,0)$ for $t\geq 0$. Consequently, the first equation of (3.16) holds. The second equation of (3.16) and the equations (3.17) are direct consequences of (3.13), (3.15) (notice that (3.13), (3.15) imply that $1=\int_{a(t)}^{b(t)}\rho(t,x)dx=(b(t)-a(t))\int_0^1\tilde{\rho}(t,\theta)d\theta$).

Therefore, for our purposes, we can consider the state to be the "transformed solution" of (2.1)-(2.6), (2.10), $(\tilde{\rho}[t],\tilde{v}[t])\in\tilde{X}$ for $t\geq 0$ given by (3.15). Notice that all equilibrium points of (2.1)-(2.6) are mapped by means of (3.15) to a single point $\left(\rho^*\chi_{[0,1]},0\chi_{[0,1]}\right)\in\tilde{X}$. This is exactly the equilibrium point that we want to stabilize globally. Stabilization is to be achieved in the topology of the linear space $X=H^1(0,1)\times C^0\left([0,1]\right)$ endowed with the norm defined by (3.14).

<u>3.4. Main results</u>

Now we are in a position to state the first main result of the present work.

**Theorem 1:** *Suppose that Assumption (H) holds. Let $\rho^*>0$ be the unique solution of equation (2.9). Then for every $R,r>0$, $S\geq 0$ there exist constants $M,\sigma>0$ (depending only on $\rho^*>0$, $R,r>0$ and $S\geq 0$) with the following property:*

**(P)** *For every classical solution of the PDE-ODE system (2.1)-(2.6) and*

$$u(t)=-R\left((r+1)v(t,a(t))+k(\rho(t,a(t)))\right), \text{ for } t>0 \quad (3.18)$$

*with $V(a(0),b(0),\rho[0],v[0])\leq S$, the following estimate holds for all $t\geq 0$:*

$$\left\|\Phi(a(t),b(t),\rho[t],v[t])-\left(\rho^*\chi_{[0,1]},0\chi_{[0,1]}\right)\right\|_X$$
$$\leq M\exp(-\sigma t)\left\|\Phi(a(0),b(0),\rho[0],v[0])-\left(\rho^*\chi_{[0,1]},0\chi_{[0,1]}\right)\right\|_X \quad (3.19)$$

Theorem 1 guarantees exponential convergence of the density profile and velocity profiles to the equilibrium (uniform) profiles in the sup-norm and the $L^2$ norm, respectively. Moreover, Theorem 1 guarantees exponential convergence to zero of the velocities of the two pistons. Estimate (3.19) shows global asymptotic stability of the equilibrium point $\left(\rho^*\chi_{[0,1]},0\chi_{[0,1]}\right)\in \tilde{X}$ (i.e., of the desired set of equilibrium points of system (2.1)-(2.6)) but not global exponential stability of the equilibrium point $\left(\rho^*\chi_{[0,1]},0\chi_{[0,1]}\right)\in \tilde{X}$ (because both the overshoot coefficient $M>0$ and the convergence rate $\sigma>0$ depend on the magnitude of $V(a(0),b(0),\rho[0],v[0])$).

Notice that since by virtue of (3.13), (3.14) we have for all $t\geq 0$

$$\left\|\Phi(a(t),b(t),\rho[t],v[t])-\left(\rho^*\chi_{[0,1]},0\chi_{[0,1]}\right)\right\|_X$$
$$=\sqrt{v^2(t,a(t))+v^2(t,b(t))}+\max_{a(t)\leq x\leq b(t)}\left(\left|\rho(t,x)-\rho^*\right|\right) \quad (3.20)$$
$$+\left(\frac{1}{b(t)-a(t)}\int_{a(t)}^{b(t)}v^2(t,x)dx\right)^{1/2}+\left((b(t)-a(t))\int_{a(t)}^{b(t)}\rho_x^2(t,x)dx\right)^{1/2}$$

it follows from (3.19) and the fact that $1=\int_{a(t)}^{b(t)}\rho(t,x)dx$ (which implies that $\frac{1}{\max_{a(t)\leq x\leq b(t)}(\rho(t,x))}\leq b(t)-a(t)\leq \frac{1}{\min_{a(t)\leq x\leq b(t)}(\rho(t,x))}$) that $\lim_{t\to+\infty}(b(t)-a(t))=(\rho^*)^{-1}$. Since (2.3), (3.19) and (3.20) imply that $\sup_{t\geq 0}(|\dot{a}(t)|\exp(\sigma t))<+\infty$ for certain $\sigma>0$ that depends on the initial conditions, it follows that there exists $a^*\in\Re$ such that $\lim_{t\to+\infty}(a(t))=a^*$.

The feedback law (3.18) in Theorem 1 requires knowledge of the function $k:(0,+\infty)\to\Re$ defined by (3.5), i.e., knowledge of the function $\mu:(0,+\infty)\to(0,+\infty)$ that describes the dependence of the viscosity of the gas on the density of the gas. Therefore, it is natural to ask the question whether a simple friction-like feedback law of the form

$$u(t)=-Rv(t,a(t)), \text{ for } t>0 \quad (3.21)$$

would achieve the same results as the feedback law (3.18). The answer is positive but requires an additional assumption for the physical properties of the gas.

**Theorem 2:** *Suppose that Assumption (H) holds. Let $\rho^*>0$ be the unique solution of equation (2.9). Moreover, suppose that the following assumption holds.*

**(A)** *There exists a constant $K>0$ such that the following inequality holds for all $\rho>0$:*

$$|k(\rho)|\leq K|P(\rho)-P(\rho^*)| \quad (3.22)$$

*Then for every $R>0$, $S\geq 0$, $r\geq RK/2$ there exist constants $M,\sigma>0$ (depending only on $\rho^*>0$, $R,K>0$ and $S\geq 0$) with the following property:*

**(P')** *For every classical solution of the PDE-ODE system (2.1)-(2.6), (3.21) with $V(a(0),b(0),\rho[0],v[0])\leq S$, estimate (3.19) holds for all $t\geq 0$.*

Assumption (A) is valid for an ideal gas under constant entropy for every $\rho^*>0$. Indeed, an ideal gas under constant entropy satisfies $P(\rho)=c\rho^\gamma$, $P=f\rho T$ and $\mu=\beta\sqrt{T}$, where $T>0$ is the temperature and $c, f, \beta, \gamma > 0$ are constants with $\gamma \in (1,2)$ (see [24]). In this case Assumption (A) is valid with any constant $K > 0$ that satisfies $K(\rho^*)^{\frac{\gamma+1}{2}} > \frac{2\beta}{(\gamma-1)\sqrt{fc}}$.

### 3.5. Auxiliary results

For the proofs of Theorem 1 and Theorem 2, we need some auxiliary lemmas. The first auxiliary lemma provides formulas for the time derivatives of the energy functional defined by (3.8), (3.9).

**Lemma 2:** *For every classical solution of the PDE-ODE system (2.1)-(2.6) the following equations hold for all $t > 0$:*

$$\frac{d}{dt} E\big(a(t), b(t), \rho[t], v[t]\big) = -\int_{a(t)}^{b(t)} \mu(\rho(t,x)) v_x^2(t,x) dx + v(t, a(t)) u(t) \tag{3.23}$$

$$\begin{aligned}
\frac{d}{dt} W\big(a(t), b(t), \rho[t], v[t]\big) = &-\int_{a(t)}^{b(t)} \rho^{-2}(t,x) P'(\rho(t,x)) \mu(\rho(t,x)) \rho_x^2(t,x) dx \\
&- k(\rho(t,b(t)))\big(P(\rho(t,b(t))) - P_{ext}\big) - k(\rho(t,a(t)))\big(P(\rho(t,a(t))) - P_{ext}\big) \\
&+ \big(v(t,a(t)) + k(\rho(t,a(t)))\big) u(t)
\end{aligned} \tag{3.24}$$

*where $E, W$ are defined by (3.8), (3.9), respectively.*

The two next auxiliary lemmas provide useful inequalities for the CLF $V$ defined by (3.7).

**Lemma 3:** *Suppose that Assumption (H) is valid. Let $S \geq 0$ and $R, r > 0$ be given. Then there exists a constant $L > 0$ (depending only on $S \geq 0$, $R, r > 0$ and $\rho^* > 0$) such that for every $a, b \in \Re$ with $a < b$, $\rho \in C^0\big([a,b];(0,+\infty)\big) \cap H^1(a,b)$, $v \in C^0\big([a,b]\big) \cap H^1(a,b)$ with $\int_a^b \rho(x) dx = 1$ and $V(a,b,\rho,v) \leq S$, where $V$ is defined by (3.7), the following inequality holds:*

$$\begin{aligned}
V(a,b,\rho,v) \leq &L\int_a^b P'(\rho(x)) \left(\frac{\partial}{\partial x}\big(k(\rho(x))\big)\right)^2 dx \\
&+ Lr\int_a^b \mu(\rho(x)) v_x^2(x) dx + LR \min\Big(\big((r+1)v(a) + k(\rho(a))\big)^2, v^2(a)\Big) \\
&+ Lk(\rho(b))\big(P(\rho(b)) - P_{ext}\big) + \frac{L}{2} k(\rho(a))\big(P(\rho(a)) - P_{ext}\big)
\end{aligned} \tag{3.25}$$

**Lemma 4:** *Suppose that Assumption (H) is valid. Let $S \geq 0$ be given. Then there exist constants $G_2 \geq G_1 > 0$ (depending only on $S \geq 0$, $r > 0$ and $\rho^* > 0$) such that for every $a, b \in \Re$ with $a < b$, $\rho \in C^0\big([a,b];(0,+\infty)\big) \cap H^1(a,b)$, $v \in C^0\big([a,b]\big)$ with $\int_a^b \rho(x) dx = 1$ and $V(a,b,\rho,v) \leq S$, where $V$ is defined by (3.7), the following inequalities hold:*

$$V(a,b,\rho,v) \leq G_2\left((b-a)\|\rho_x\|_2^2 + \|\rho - \rho^*\chi_{[a,b]}\|_\infty^2 + v^2(a) + v^2(b) + (b-a)^{-1}\|v\|_2^2\right) \quad (3.26)$$

$$G_1\left((b-a)\|\rho_x\|_2^2 + \|\rho - \rho^*\chi_{[a,b]}\|_\infty^2 + v^2(a) + v^2(b) + (b-a)^{-1}\|v\|_2^2\right) \leq V(a,b,\rho,v)$$

where $\|\rho - \rho^*\chi_{[a,b]}\|_\infty = \max_{a \leq x \leq b}(|\rho(x) - \rho^*|)$, $\|v\|_2 = \left(\int_a^b v^2(x)dx\right)^{1/2}$ and $\|\rho_x\|_2 = \left(\int_a^b \rho_x^2(x)dx\right)^{1/2}$.

Inequality (3.25) provides an estimate of the dissipation rate of the Lyapunov functional for the closed-loop systems (2.1)-(2.6) with (3.18) or (3.21). On the other hand, inequalities (3.26) provide estimates of the Lyapunov functional in terms of the norm of the state space.

## 4. Proofs of Main Results

We start by providing the proof of Lemma 1.

**Proof of Lemma 1:** Let $S \geq 0$ be given (arbitrary) and let $a, b \in \mathfrak{R}$ with $a < b$, $\rho \in C^0([a,b];(0,+\infty)) \cap H^1(a,b)$, $v \in C^0([a,b])$ with $\int_a^b \rho(x)dx = 1$ be given (arbitrary). Suppose that (3.11) holds.

Since $V(a,b,\rho,v) \leq S$ we have (from (3.7), (3.9)) that $\int_a^b \rho(x)\left(v(x) + \rho^{-2}(x)\mu(\rho(x))\rho_x(x)\right)^2 dx \leq 2S$.

Moreover, since $(x+y)^2 \geq \frac{1}{2}y^2 - x^2$ for all $x, y \in \mathfrak{R}$, it holds that

$$\int_a^b \rho^{-3}(x)\mu^2(\rho(x))\rho_x^2(x)dx \leq 4S + 2\int_a^b \rho(x)v^2(x)dx \quad (4.1)$$

Furthermore, since $V(a,b,\rho,v) \leq S$, it follows (from (3.7), (3.8)) that $\int_a^b \rho(x)v^2(x)dx \leq 2r^{-1}S$. Consequently, combining (4.1) with the previous inequality we get

$$\int_a^b \rho^{-3}(x)\mu^2(\rho(x))\rho_x^2(x)dx \leq 4(1+r^{-1})S \quad (4.2)$$

The function $G$ defined by (3.4) satisfies the following inequalities for all $x, y \in [a,b]$ (by virtue of the Cauchy-Schwarz inequality):

$$|G(\rho(x)) - G(\rho(y))| = \left|\int_y^x G'(\rho(s))\rho_s(s)ds\right|$$

$$\leq \left(\int_y^x \rho^{-3}(s)\mu^2(\rho(s))\rho_s^2(s)ds\right)^{1/2}\left(\int_y^x \frac{\rho^3(s)(G'(\rho(s)))^2}{\mu^2(\rho(s))}ds\right)^{1/2} \quad (4.3)$$

$$\leq \left(\int_a^b \rho^{-3}(s)\mu^2(\rho(s))\rho_s^2(s)ds\right)^{1/2}\left(\int_a^b Q(\rho(s))ds\right)^{1/2}$$

Since $V(a,b,\rho,v) \leq S$, it follows (from (3.7), (3.8), (3.9), (3.10)) that $\int_a^b Q(\rho(x))dx \leq \dfrac{S}{r+1}$. Combining the previous inequality with (4.2), (4.3) gives:

$$|G(\rho(x)) - G(\rho(y))| \leq \frac{2}{\sqrt{r}} S, \text{ for all } x, y \in [a,b] \tag{4.4}$$

Since $V(a,b,\rho,v) \leq S$, it follows (from (3.7), (3.9)) that $(v(a) + k(\rho(a)))^2 \leq 2S$. Moreover, since $(x+y)^2 \geq \dfrac{1}{2} y^2 - x^2$ for all $x, y \in \Re$, it holds that $k^2(\rho(a)) \leq 4S + 2v^2(a)$. Since $V(a,b,\rho,v) \leq S$, it follows (from (3.7), (3.8)) that $v^2(a) \leq 2r^{-1}S$. Combining the previous inequalities we obtain the following estimate:

$$k^2(\rho(a)) \leq 4(1 + r^{-1})S \tag{4.5}$$

Equation (3.3) and estimate (4.5) imply the existence of $\rho_1 > 0$ (that depends only on $S \geq 0$, $r > 0$ and $\rho^* > 0$) such that $\rho(a) \leq \rho_1$. Using (4.4) with $y = a$ and the fact that $G$ is increasing, we get for all $x \in [a,b]$:

$$G(\rho(x)) \leq G(\rho(a)) + 2S \leq G(\rho_1) + \frac{2}{\sqrt{r}} S \tag{4.6}$$

Inequality (4.6) in conjunction with (3.1) shows the existence of $\rho_{\max} > 0$ (that depends only on $S \geq 0$, $r > 0$ and $\rho^* > 0$) such that the right inequality (3.12) holds.

Using the fact that $\int_a^b \rho(x)dx = 1$ we have:

$$\min_{x \in [a,b]} (\rho(x)) \leq \frac{1}{b-a} \leq \max_{x \in [a,b]} (\rho(x)) \tag{4.7}$$

We next consider the following cases:

<u>Case 1:</u> $\rho^*/2 < \min_{x \in [a,b]} (\rho(x)) \leq \max_{x \in [a,b]} (\rho(x))$. In this case the inequality $\min_{x \in [a,b]} (\rho(x)) \geq \rho^*/2$ holds.

<u>Case 2:</u> $\min_{x \in [a,b]} (\rho(x)) \leq \rho^*/2 \leq \max_{x \in [a,b]} (\rho(x))$. By continuity there exists $x^* \in [a,b]$ such that $\rho(x^*) = \rho^*/2$. Applying inequality (4.4) with $y = x^*$ we get for all $x \in [a,b]$:

$$G(\rho(x)) \geq G(\rho(x^*)) - \frac{2}{\sqrt{r}} S = G\left(\frac{\rho^*}{2}\right) - \frac{2}{\sqrt{r}} S \tag{4.8}$$

By virtue of (3.2) there exists $\rho_2 \in \left(0, \dfrac{\rho^*}{2}\right]$ such that $G\left(\dfrac{\rho^*}{2}\right) - \dfrac{2}{\sqrt{r}} S = G(\rho_2)$. Inequality (4.8) and the fact that $G$ is increasing, imply that $\min_{x \in [a,b]} (\rho(x)) \geq \rho_2$.

<u>Case 3:</u> $\min_{x \in [a,b]} (\rho(x)) \leq \max_{x \in [a,b]} (\rho(x)) < \rho^*/2$. By continuity and (4.7) there exists $x^* \in [a,b]$ such that $\rho(x^*) = \dfrac{1}{b-a}$. Applying inequality (4.4) with $y = x^*$, we get for all $x \in [a,b]$

$$G(\rho(x)) \geq G(\rho(x^*)) - \frac{2}{\sqrt{r}} S = G\left(\frac{1}{b-a}\right) - \frac{2}{\sqrt{r}} S \quad (4.9)$$

Since $V(a,b,\rho,v) \leq S$, it follows (from (3.7), (3.8), (3.9), (3.10)) that $(b-a)\min_{x \in [a,b]}(Q(\rho(x))) \leq \int_a^b Q(\rho(x))dx \leq \frac{S}{r+1}$. Consequently, in this case we obtain the following estimate

$$\frac{r+1}{S+1} \min_{x \in [a,b]}(Q(\rho(x))) \leq \frac{1}{b-a} \quad (4.10)$$

Since $\min_{x \in [a,b]}(\rho(x)) \leq \max_{x \in [a,b]}(\rho(x)) \leq \frac{\rho^*}{2}$ and since $Q''(\rho) = \frac{P'(\rho)}{\rho} > 0$ for all $\rho > 0$, $Q(\rho^*) = Q'(\rho^*) = 0$ (recall (3.6)), it follows that $\min_{x \in [a,b]}(Q(\rho(x))) \geq Q\left(\frac{\rho^*}{2}\right)$. Therefore, we get from (4.10) the inequality $\frac{r+1}{S+1} Q\left(\frac{\rho^*}{2}\right) \leq \frac{1}{b-a}$. Since $G$ is increasing we have $G\left(\frac{r+1}{S+1} Q\left(\frac{\rho^*}{2}\right)\right) \leq G\left(\frac{1}{b-a}\right)$. Consequently, we get from (4.9) for all $x \in [a,b]$:

$$G(\rho(x)) \geq G\left(\frac{r+1}{S+1} Q\left(\frac{\rho^*}{2}\right)\right) - \frac{2}{\sqrt{r}} S \quad (4.11)$$

By virtue of (3.1), (3.2) there exists $\rho_3 \in (0, +\infty)$ such that $G\left(\frac{r+1}{S+1} Q\left(\frac{\rho^*}{2}\right)\right) - \frac{2}{\sqrt{r}} S = G(\rho_3)$. Inequality (4.11) and the fact that $G$ is increasing, imply that $\min_{x \in [a,b]}(\rho(x)) \geq \rho_3$.

Combining all cases, we conclude that the left inequality (3.12) holds with $\rho_{\min} = \min(\rho_2, \rho_3)$. The proof is complete. ◁

We continue with the proof of Lemma 2.

**Proof of Lemma 2:** Consider a classical solution for the PDE-ODE system (2.1)-(2.6), i.e., we consider functions $a, b \in C^1(\Re_+) \cap C^2((0, +\infty))$ satisfying $b(t) > a(t)$ for all $t \geq 0$, $\rho \in C^1(\bar{\Omega}; (0, +\infty)) \cap C^2(\Omega)$, $v \in C^0(\bar{\Omega}) \cap C^1(\tilde{\Omega})$ with $v[t] \in C^2((a(t), b(t)))$ for each $t > 0$ and $\Omega = \bigcup_{t>0}(\{t\} \times (a(t), b(t)))$, $\tilde{\Omega} = \bigcup_{t>0}(\{t\} \times [a(t), b(t)])$ that satisfy equations (2.1)-(2.6) for a given input $u \in C^0(\Re_+)$. Using (2.3), (3.8) and (3.10) we obtain for all $t \geq 0$:

$$E(a(t), b(t), \rho[t], v[t]) = \frac{\dot{a}^2(t)}{2} + \frac{\dot{b}^2(t)}{2} + \frac{1}{2} \int_{a(t)}^{b(t)} \rho(t,x) v^2(t,x) dx + \int_{a(t)}^{b(t)} Q(\rho(t,x)) dx \quad (4.12)$$

It follows from (4.12) that the following equation holds for all $t > 0$:

$$\frac{d}{dt}E(a(t),b(t),\rho[t],v[t]) = \dot{a}(t)\ddot{a}(t) + \dot{b}(t)\ddot{b}(t) + \frac{1}{2}\int_{a(t)}^{b(t)} \rho_t(t,x)v^2(t,x)dx$$

$$+ \frac{1}{2}\rho(t,b(t))v^2(t,b(t))\dot{b}(t) - \rho(t,a(t))v^2(t,a(t))\dot{a}(t) + \int_{a(t)}^{b(t)} Q'(\rho(t,x))\rho_t(t,x)dx \quad (4.13)$$

$$+ \int_{a(t)}^{b(t)} \rho(t,x)v(t,x)v_t(t,x)dx + Q(\rho(t,b(t)))\dot{b}(t) - Q(\rho(t,a(t)))\dot{a}(t)$$

Using (2.1), (2.3), (2.4) and (2.5) we get from (4.13) for all $t > 0$:

$$\frac{d}{dt}E(a(t),b(t),\rho[t],v[t]) = \dot{a}(t)\big(u(t) + P_{ext} - P(\rho(t,a(t))) + \mu(\rho(t,a(t)))v_x(t,a(t))\big)$$

$$+ \dot{b}(t)\big(P(\rho(t,b(t))) - \mu(\rho(t,b(t)))v_x(t,b(t)) - P_{ext}\big) - \int_{a(t)}^{b(t)} Q'(\rho(t,x))(\rho v)_x(t,x)dx$$

$$- \frac{1}{2}\int_{a(t)}^{b(t)} (\rho v)_x(t,x)v^2(t,x)dx + \frac{1}{2}\rho(t,b(t))\dot{b}^3(t) - \rho(t,a(t))\dot{a}^3(t) \quad (4.14)$$

$$+ \int_{a(t)}^{b(t)} \rho(t,x)v(t,x)v_t(t,x)dx + Q(\rho(t,b(t)))\dot{b}(t) - Q(\rho(t,a(t)))\dot{a}(t)$$

Integrating by parts and using (2.3), we obtain from (4.14) for all $t > 0$:

$$\frac{d}{dt}E(a(t),b(t),\rho[t],v[t]) = \dot{a}(t)\big(u(t) + P_{ext} - P(\rho(t,a(t))) + \mu(\rho(t,a(t)))v_x(t,a(t))\big)$$

$$+ \dot{b}(t)\big(P(\rho(t,b(t))) - \mu(\rho(t,b(t)))v_x(t,b(t)) - P_{ext}\big)$$

$$+ \int_{a(t)}^{b(t)} \rho(t,x)v(t,x)\big(v_t(t,x) + v(t,x)v_x(t,x) + Q''(\rho(t,x))\rho_x(t,x)\big)dx \quad (4.15)$$

$$+ \big(Q(\rho(t,b(t))) - Q'(\rho(t,b(t)))\rho(t,b(t))\big)\dot{b}(t)$$

$$+ \big(Q'(\rho(t,a(t)))\rho(t,a(t)) - Q(\rho(t,a(t)))\big)\dot{a}(t)$$

Definition (3.6) implies the equations $Q(\rho) - \rho Q'(\rho) = P(\rho^*) - P(\rho)$, $Q''(\rho) = P'(\rho)/\rho$ for all $\rho > 0$. The previous equations in conjunction with (4.15), (2.2) and (2.9) give the following equation for all $t > 0$:

$$\frac{d}{dt}E(a(t),b(t),\rho[t],v[t]) = \dot{a}(t)\big(u(t) + \mu(\rho(t,a(t)))v_x(t,a(t))\big)$$

$$- \mu(\rho(t,b(t)))v_x(t,b(t))\dot{b}(t) + \int_{a(t)}^{b(t)} v(t,x)\frac{\partial}{\partial x}\big(\mu(\rho(t,x))v_x(t,x)\big)dx \quad (4.16)$$

Integrating by parts and using (4.16) and (2.3), we get (3.23).

Define for all $(t,x) \in \overline{\Omega} = \bigcup_{t \geq 0}\big(\{t\} \times [a(t),b(t)]\big)$:

$$w(t,x) := \rho(t,x)v(t,x) + k'(\rho(t,x))\rho_x(t,x) \quad (4.17)$$

Using (2.1), (2.2) and definition (3.5), we obtain from (4.17) for all $(t,x) \in \Omega = \bigcup_{t > 0}\big(\{t\} \times (a(t),b(t))\big)$:

$$w_t(t,x) = -v_x(t,x)w(t,x) - w_x(t,x)v(t,x) - P'(\rho(t,x))\rho_x(t,x) \quad (4.18)$$

Using (3.9), (3.10), (2.3) and (4.17), we obtain for all $t \geq 0$:

$$W(a(t),b(t),\rho[t],v[t]) = \frac{1}{2}\int_{a(t)}^{b(t)} \rho^{-1}(t,x)w^2(t,x)dx$$
$$+ \frac{1}{2}\left(\dot{b}(t) - k(\rho(t,b(t)))\right)^2 + \frac{1}{2}\left(\dot{a}(t) + k(\rho(t,a(t)))\right)^2 + \int_{a(t)}^{b(t)} Q(\rho(t,x))dx \quad (4.19)$$

Differentiating with respect to time we obtain from (4.19) for all $t > 0$:

$$\frac{d}{dt}W(a(t),b(t),\rho[t],v[t]) = \frac{1}{2}\rho^{-1}(t,b(t))w^2(t,b(t))\dot{b}(t)$$
$$- \frac{1}{2}\rho^{-1}(t,a(t))w^2(t,a(t))\dot{a}(t)$$
$$- \frac{1}{2}\int_{a(t)}^{b(t)} \rho^{-2}(t,x)\rho_t(t,x)w^2(t,x)dx + \int_{a(t)}^{b(t)} \rho^{-1}(t,x)w(t,x)w_t(t,x)dx$$
$$+ \left(\dot{b}(t) - k(\rho(t,b(t)))\right)\frac{d}{dt}\left(\dot{b}(t) - k(\rho(t,b(t)))\right) \quad (4.20)$$
$$+ \left(\dot{a}(t) + k(\rho(t,a(t)))\right)\frac{d}{dt}\left(\dot{a}(t) + k(\rho(t,a(t)))\right)$$
$$+ Q(\rho(t,b(t)))\dot{b}(t) - Q(\rho(t,b(t)))\dot{a}(t) + \int_{a(t)}^{b(t)} Q'(\rho(t,x))\rho_t(t,x)dx$$

Using (2.1), (2.3), (2.4), (2.5) and definition (3.5), we get for all $t > 0$:

$$\frac{d}{dt}\left(\dot{b}(t) - k(\rho(t,b(t)))\right) = \ddot{b}(t) - k'(\rho(t,b(t)))\frac{d}{dt}\left(\rho(t,b(t))\right)$$
$$= \ddot{b}(t) - \frac{\mu(\rho(t,b(t)))}{\rho(t,b(t))}\left(\frac{\partial \rho}{\partial t}(t,b(t)) + \frac{\partial \rho}{\partial x}(t,b(t))\dot{b}(t)\right) \quad (4.21)$$
$$= \ddot{b}(t) + \mu(\rho(t,b(t)))v_x(t,b(t))$$
$$= P(\rho(t,b(t))) - P_{ext}$$

$$\frac{d}{dt}\left(\dot{a}(t) + k(\rho(t,a(t)))\right) = \ddot{a}(t) + k'(\rho(t,a(t)))\frac{d}{dt}\left(\rho(t,a(t))\right)$$
$$= \ddot{a}(t) + \frac{\mu(\rho(t,a(t)))}{\rho(t,a(t))}\left(\frac{\partial \rho}{\partial t}(t,a(t)) + \frac{\partial \rho}{\partial x}(t,a(t))\dot{a}(t)\right) \quad (4.22)$$
$$= \ddot{a}(t) - \mu(\rho(t,a(t)))v_x(t,a(t))$$
$$= u(t) + P_{ext} - P(\rho(t,a(t)))$$

Combining (4.20), (4.21), (4.22), (2.1) and (4.18) we get for all $t > 0$:

$$\frac{d}{dt}W(a(t),b(t),\rho[t],v[t]) = \frac{1}{2}\rho^{-1}(t,b(t))w^2(t,b(t))\dot{b}(t)$$

$$-\frac{1}{2}\rho^{-1}(t,a(t))w^2(t,a(t))\dot{a}(t)$$

$$+\frac{1}{2}\int_{a(t)}^{b(t)}\rho^{-2}(t,x)(\rho v)_x(t,x)w^2(t,x)dx - \int_{a(t)}^{b(t)}\rho^{-1}(t,x)w(t,x)(vw)_x(t,x)dx$$

$$-\int_{a(t)}^{b(t)}\rho^{-1}(t,x)w(t,x)P'(\rho(t,x))\rho_x(t,x)dx \qquad (4.23)$$

$$+\left(\dot{b}(t)-k(\rho(t,b(t)))\right)\left(P(\rho(t,b(t)))-P_{ext}\right)$$

$$+\left(\dot{a}(t)+k(\rho(t,a(t)))\right)\left(u(t)+P_{ext}-P(\rho(t,a(t)))\right)$$

$$+Q(\rho(t,b(t)))\dot{b}(t)-Q(\rho(t,b(t)))\dot{a}(t)-\int_{a(t)}^{b(t)}Q'(\rho(t,x))(\rho v)_x(t,x)dx$$

Integrating by parts and using (2.3), we obtain from (4.23) for all $t > 0$:

$$\frac{d}{dt}W(a(t),b(t),\rho[t],v[t]) = \rho^{-1}(t,b(t))w^2(t,b(t))\dot{b}(t)$$

$$-\rho^{-1}(t,a(t))w^2(t,a(t))\dot{a}(t)$$

$$+\int_{a(t)}^{b(t)}\rho^{-2}(t,x)v(t,x)\rho_x(t,x)w^2(t,x)dx - 2\int_{a(t)}^{b(t)}\rho^{-1}(t,x)v(t,x)w(t,x)w_x(t,x)dx$$

$$-\int_{a(t)}^{b(t)}\rho^{-1}(t,x)w^2(t,x)v_x(t,x)dx - \int_{a(t)}^{b(t)}\rho^{-1}(t,x)w(t,x)P'(\rho(t,x))\rho_x(t,x)dx \qquad (4.24)$$

$$+\left(\dot{b}(t)-k(\rho(t,b(t)))\right)\left(P(\rho(t,b(t)))-P_{ext}\right)$$

$$+\left(\dot{a}(t)+k(\rho(t,a(t)))\right)\left(u(t)+P_{ext}-P(\rho(t,a(t)))\right)$$

$$+Q(\rho(t,b(t)))\dot{b}(t)-Q(\rho(t,b(t)))\dot{a}(t)-Q'(\rho(t,b(t)))\rho(t,b(t))v(t,b(t))$$

$$+Q'(\rho(t,a(t)))\rho(t,a(t))v(t,a(t))+\int_{a(t)}^{b(t)}Q''(\rho(t,x))\rho_x(t,x)\rho(t,x)v(t,x)dx$$

Definition (3.6) implies the equation $Q''(\rho) = P'(\rho)/\rho$ for all $\rho > 0$. The previous equation in conjunction with (4.24) gives the following equation for all $t > 0$:

$$\frac{d}{dt}W(a(t),b(t),\rho[t],v[t]) = \rho^{-1}(t,b(t))w^2(t,b(t))\dot{b}(t)$$

$$-\rho^{-1}(t,a(t))w^2(t,a(t))\dot{a}(t)$$

$$-\int_{a(t)}^{b(t)}\frac{\partial}{\partial x}\left(\rho^{-1}(t,x)w^2(t,x)v(t,x)\right)dx$$

$$+\left(\dot{b}(t)-k(\rho(t,b(t)))\right)\left(P(\rho(t,b(t)))-P_{ext}\right) \qquad (4.25)$$

$$+\left(\dot{a}(t)+k(\rho(t,a(t)))\right)\left(u(t)+P_{ext}-P(\rho(t,a(t)))\right)$$

$$+Q(\rho(t,b(t)))\dot{b}(t)-Q(\rho(t,b(t)))\dot{a}(t)-Q'(\rho(t,b(t)))\rho(t,b(t))\dot{b}(t)$$

$$+Q'(\rho(t,a(t)))\rho(t,a(t))\dot{a}(t)+\int_{a(t)}^{b(t)}P'(\rho(t,x))\rho_x(t,x)\left(v(t,x)-\rho^{-1}(t,x)w(t,x)\right)dx$$

Definition (3.6) implies the equation $Q(\rho) - \rho Q'(\rho) = P(\rho^*) - P(\rho)$ for all $\rho > 0$. The previous equation in conjunction with (4.25) and (2.9) give the following equation for all $t > 0$:

$$\frac{d}{dt}W(a(t), b(t), \rho[t], v[t]) = -k(\rho(t,b(t)))\left(P(\rho(t,b(t))) - P_{ext}\right)$$
$$-k(\rho(t,a(t)))\left(P(\rho(t,a(t))) - P_{ext}\right) + \left(\dot{a}(t) + k(\rho(t,a(t)))\right)u(t) \quad (4.26)$$
$$+ \int_{a(t)}^{b(t)} P'(\rho(t,x))\rho_x(t,x)\left(v(t,x) - \rho^{-1}(t,x)w(t,x)\right)dx$$

Equation (3.24) is a direct consequence of (4.26) and definitions (3.5), (4.17).
The proof is complete. ◁

Next, we provide the proofs of Lemma 3 and Lemma 4.

**Proof of Lemma 3:** Let $S \geq 0$ be given (arbitrary) and let $a, b \in \Re$ with $a < b$, $\rho \in C^0([a,b]; (0,+\infty)) \cap H^1(a,b)$, $v \in C^0([a,b])$ with $\int_a^b \rho(x)dx = 1$ be given (arbitrary). Suppose that (3.11) holds. By virtue of Lemma 1, there exist constants $\rho_{max} \geq \rho_{min} > 0$ (depending only on $S \geq 0$, $r > 0$ and $\rho^* > 0$) such that (3.12) holds.

Definition (3.5) and the fact that $P'(\rho) > 0$ for all $\rho > 0$ imply that there exists a constant $K > 0$ (depending only on $S \geq 0$, $r > 0$ and $\rho^* > 0$) such that the following inequality holds for all $\rho \in [\rho_{min}, \rho_{max}]$:

$$|k(\rho)| \leq K|P(\rho) - P(\rho^*)| \quad (4.27)$$

Using the fact that $k(\rho)(P(\rho) - P(\rho^*)) \geq 0$ for all $\rho > 0$ (recall (3.5) and the fact that $P$ is increasing) in conjunction with (4.27) and the inequality $2v(a)k(\rho(a)) \geq -\varepsilon^{-1}v^2(a) - \varepsilon|k(\rho(a))|^2$ with $\varepsilon = \dfrac{1+4KR}{4KR(r+1)}$, we obtain from (2.9) the following inequalities:

$$R((r+1)v(a) + k(\rho(a)))^2 + k(\rho(b))(P(\rho(b)) - P_{ext}) + k(\rho(a))(P(\rho(a)) - P_{ext})$$
$$\geq \frac{1}{K}|k(\rho(a))|^2 + \frac{1}{K}|k(\rho(b))|^2 + R((r+1)v(a) + k(\rho(a)))^2$$
$$= \left(\frac{1}{K} + R\right)|k(\rho(a))|^2 + \frac{1}{K}|k(\rho(b))|^2 + R(r+1)^2 v^2(a) + 2R(r+1)v(a)k(\rho(a)) \quad (4.28)$$
$$\geq (R + \frac{1}{K} - \varepsilon R(r+1))|k(\rho(a))|^2 + \frac{1}{K}|k(\rho(b))|^2 + \varepsilon^{-1}R(r+1)(\varepsilon(r+1) - 1)v^2(a)$$
$$\geq \frac{3}{4K}|k(\rho(a))|^2 + \frac{1}{K}|k(\rho(b))|^2 + \frac{R(r+1)^2}{1+4KR}v^2(a)$$

$$Rv^2(a) + k(\rho(b))(P(\rho(b)) - P_{ext}) + k(\rho(a))(P(\rho(a)) - P_{ext})$$
$$\geq \frac{1}{K}|k(\rho(a))|^2 + \frac{1}{K}|k(\rho(b))|^2 + Rv^2(a) \tag{4.29}$$

Inequalities (4.28), (4.29) allow us to conclude that there exists a constant $c_1 > 0$ (depending only on $S \geq 0$, $r, R > 0$ and $\rho^* > 0$) such that the following inequality holds:

$$R \min\left(((r+1)v(a) + k(\rho(a)))^2, v^2(a)\right)$$
$$+ k(\rho(b))(P(\rho(b)) - P_{ext}) + \frac{1}{2}k(\rho(a))(P(\rho(a)) - P_{ext}) \tag{4.30}$$
$$\geq c_1\left(|k(\rho(a))|^2 + |k(\rho(b))|^2 + v^2(a)\right)$$

Using the inequality $(x+y)^2 \leq 2x^2 + 2y^2$ that holds for all $x, y \in \Re$ and definitions (3.7), (3.8), (3.9), (3.10), we obtain the following inequality:

$$V(a,b,\rho,v) \leq \int_a^b \rho^{-3}(x)\mu^2(\rho(x))\rho_x^2(x)dx$$
$$+ |k(\rho(b))|^2 + |k(\rho(a))|^2 + (r+1)\int_a^b Q(\rho(x))dx \tag{4.31}$$
$$+ (r+2)\frac{v^2(a)}{2} + (r+2)\frac{v^2(b)}{2} + \frac{r+2}{2}\int_a^b \rho(x)v^2(x)dx$$

The Cauchy-Schwarz inequality and the triangle inequality allows us to conclude that the following implications are true for all $x \in [a,b]$:

$$v(x) = v(a) + \int_a^x v_s(s)ds$$
$$\Rightarrow |v(x)| \leq |v(a)| + \left(\int_a^x \frac{ds}{\mu(\rho(s))}\right)^{1/2}\left(\int_a^x \mu(\rho(s))v_s^2(s)ds\right)^{1/2}$$
$$\Rightarrow |v(x)| \leq |v(a)| + \left(\int_a^b \frac{ds}{\mu(\rho(s))}\right)^{1/2}\left(\int_a^b \mu(\rho(s))v_s^2(s)ds\right)^{1/2}$$

Using the above inequality in conjunction with (3.12) and the inequality $(x+y)^2 \leq 2x^2 + 2y^2$ that holds for all $x, y \in \Re$, we get for all $x \in [a,b]$:

$$v^2(x) \leq 2v^2(a) + 2\Gamma(b-a)\int_a^b \mu(\rho(s))v_s^2(s)ds \tag{4.32}$$

where $\Gamma := \dfrac{1}{\min\{\mu(\rho): \rho_{\min} \leq \rho \leq \rho_{\max}\}}$ (a constant that depends only on $S \geq 0$, $r > 0$ and $\rho^* > 0$). Using the fact that $\int_a^b \rho(x)dx = 1$ in conjunction with (3.12) gives:

$$\frac{1}{\rho_{\max}} \leq b - a \leq \frac{1}{\rho_{\min}} \tag{4.33}$$

Therefore, it follows from (4.32) and (4.33) that the following inequality holds:

$$\max\{v^2(x): x \in [a,b]\} \leq 2v^2(a) + \frac{2\Gamma}{\rho_{\min}} \int_a^b \mu(\rho(x)) v_x^2(x) dx \tag{4.34}$$

Combining (4.31) with (3.12), (4.33) and (4.34) we get:

$$\begin{aligned}
V(a,b,\rho,v) &\leq \int_a^b \rho^{-3}(x) \mu^2(\rho(x)) \rho_x^2(x) dx \\
&+ |k(\rho(b))|^2 + |k(\rho(a))|^2 + (r+1) \int_a^b Q(\rho(x)) dx \\
&+ \frac{\rho_{\max}}{\rho_{\min}}(r+4)^2 v^2(a) + (r+2) \frac{2\Gamma \rho_{\max}}{\rho_{\min}^2} \int_a^b \mu(\rho(x)) v_x^2(x) dx
\end{aligned} \tag{4.35}$$

Using (4.35) in conjunction with (3.12), we get:

$$\begin{aligned}
V(a,b,\rho,v) &\leq c_2 \int_a^b \rho^{-2}(x) P'(\rho(x)) \mu(\rho(x)) \rho_x^2(x) dx \\
&+ |k(\rho(b))|^2 + |k(\rho(a))|^2 + (r+1) \int_a^b Q(\rho(x)) dx \\
&+ \frac{\rho_{\max}}{\rho_{\min}}(r+4)^2 v^2(a) + (r+2) \frac{2\Gamma \rho_{\max}}{\rho_{\min}^2} \int_a^b \mu(\rho(x)) v_x^2(x) dx
\end{aligned} \tag{4.36}$$

where $c_2 := \max\left\{\frac{\mu(\tau)}{\tau P'(\tau)} : \rho_{\min} \leq \tau \leq \rho_{\max}\right\}$ (a constant that depends only on $S \geq 0$, $r > 0$ and $\rho^* > 0$).

The Cauchy-Schwarz inequality allows us to conclude that the following inequalities hold for all $x \in [a,b]$:

$$\begin{aligned}
Q(\rho(x)) &= Q(\rho(a)) + \int_a^x Q'(\rho(s)) \rho_s(s) ds \\
&\leq Q(\rho(a)) + \left(\int_a^x \frac{\rho^2(s)(Q'(\rho(s)))^2}{P'(\rho(s))\mu(\rho(s))} ds\right)^{1/2} \left(\int_a^x \rho^{-2}(s) P'(\rho(s))\mu(\rho(s)) \rho_s^2(s) ds\right)^{1/2} \\
&\leq Q(\rho(a)) + \left(\int_a^b \frac{\rho^2(s)(Q'(\rho(s)))^2}{P'(\rho(s))\mu(\rho(s))} ds\right)^{1/2} \left(\int_a^b \rho^{-2}(s) P'(\rho(s))\mu(\rho(s)) \rho_s^2(s) ds\right)^{1/2}
\end{aligned} \tag{4.37}$$

Using (4.37), (4.33) and the inequality $xy \leq \frac{\varepsilon}{2} x^2 + \frac{1}{2\varepsilon} y^2$ that holds for all $x, y \in \Re$ and for all $\varepsilon > 0$, we obtain the following estimates for all $\varepsilon > 0$:

$$\int_a^b Q(\rho(x))dx \leq (b-a)Q(\rho(a))$$

$$+(b-a)\left(\int_a^b \frac{\rho^2(x)(Q'(\rho(x)))^2}{P'(\rho(x))\mu(\rho(x))}dx\right)^{1/2}\left(\int_a^b \rho^{-2}(x)P'(\rho(x))\mu(\rho(x))\rho_x^2(x)dx\right)^{1/2}$$

$$\leq \rho_{\min}^{-1}Q(\rho(a))$$

$$+\rho_{\min}^{-1}\left(\int_a^b \frac{\rho^2(x)(Q'(\rho(x)))^2}{P'(\rho(x))\mu(\rho(x))}dx\right)^{1/2}\left(\int_a^b \rho^{-2}(x)P'(\rho(x))\mu(\rho(x))\rho_x^2(x)dx\right)^{1/2} \quad (4.38)$$

$$\leq \rho_{\min}^{-1}Q(\rho(a)) + \frac{\varepsilon}{2}\int_a^b \frac{\rho^2(x)(Q'(\rho(x)))^2}{P'(\rho(x))\mu(\rho(x))}dx$$

$$+\frac{1}{2\varepsilon\rho_{\min}^2}\int_a^b \rho^{-2}(x)P'(\rho(x))\mu(\rho(x))\rho_x^2(x)dx$$

Definition (3.6) implies that $Q''(\rho) = \rho^{-1}P'(\rho) > 0$ for all $\rho > 0$, $Q(\rho^*) = Q'(\rho^*) = 0$. Consequently, it follows that $Q(\rho) \geq \frac{1}{2}(\rho - \rho^*)^2 \min\{s^{-1}P'(s) : \rho_{\min} \leq s \leq \rho_{\max}\}$ and $|Q'(\rho)| \leq |\rho - \rho^*|\max\{s^{-1}P'(s) : \rho_{\min} \leq s \leq \rho_{\max}\}$ for all $\rho \in [\rho_{\min}, \rho_{\max}]$. Therefore, there exists a constant $\tilde{K} > 0$ (depending only on $S \geq 0$, $r > 0$ and $\rho^* > 0$) such that:

$$\frac{\rho^2(Q'(\rho))^2}{P'(\rho)\mu(\rho)} \leq \tilde{K}Q(\rho), \text{ for all } \rho \in [\rho_{\min}, \rho_{\max}] \quad (4.39)$$

Using (4.38), (4.39) and (3.12), we conclude that the following estimate holds for all $\varepsilon > 0$:

$$\int_a^b Q(\rho(x))dx \leq \rho_{\min}^{-1}Q(\rho(a)) + \frac{\varepsilon\tilde{K}}{2}\int_a^b Q(\rho(x))dx$$

$$+\frac{1}{2\varepsilon\rho_{\min}^2}\int_a^b \rho^{-2}(x)P'(\rho(x))\mu(\rho(x))\rho_x^2(x)dx \quad (4.40)$$

Selecting $\varepsilon = 1/\tilde{K}$, we get from (4.40):

$$\int_a^b Q(\rho(x))dx \leq \frac{2}{\rho_{\min}}Q(\rho(a)) + \frac{\tilde{K}}{\rho_{\min}^2}\int_a^b \rho^{-2}(x)P'(\rho(x))\mu(\rho(x))\rho_x^2(x)dx \quad (4.41)$$

Combining (4.36) and (4.41), we obtain:

$$V(a,b,\rho,v) \leq \left(c_2 + \frac{\tilde{K}(r+1)}{\rho_{\min}^2}\right)\int_a^b \rho^{-2}(x)P'(\rho(x))\mu(\rho(x))\rho_x^2(x)dx$$

$$+|k(\rho(b))|^2 + |k(\rho(a))|^2 + \frac{2(r+1)}{\rho_{\min}}Q(\rho(a)) \quad (4.42)$$

$$+\frac{\rho_{\max}}{\rho_{\min}}(r+4)^2 v^2(a) + (r+2)\frac{2\Gamma\rho_{\max}}{\rho_{\min}^2}\int_a^b \mu(\rho(x))v_x^2(x)dx$$

Definitions (3.5), (3.6) imply that $Q''(\rho) = \rho^{-1} P'(\rho) > 0$, $k'(\rho) = \rho^{-1} \mu(\rho) > 0$ for all $\rho > 0$ and $Q(\rho^*) = Q'(\rho^*) = k(\rho^*) = 0$. Consequently, it follows that $Q(\rho) \leq \frac{1}{2}(\rho - \rho^*)^2 \max\{s^{-1} P'(s): \rho_{\min} \leq s \leq \rho_{\max}\}$ and $|k(\rho)| \geq |\rho - \rho^*| \min\{s^{-1} \mu(s): \rho_{\min} \leq s \leq \rho_{\max}\}$ for all $\rho \in [\rho_{\min}, \rho_{\max}]$. Therefore, there exists a constant $F > 0$ (depending only on $S \geq 0$, $r > 0$ and $\rho^* > 0$) such that:

$$Q(\rho) \leq F |k(\rho)|^2, \text{ for all } \rho \in [\rho_{\min}, \rho_{\max}] \quad (4.43)$$

Using (4.42), (4.43) and (3.12), we conclude that the following estimate holds:

$$V(a,b,\rho,v) \leq \left(c_2 + \frac{\tilde{K}(r+1)}{\rho_{\min}^2}\right) \int_a^b \rho^{-2}(x) P'(\rho(x)) \mu(\rho(x)) \rho_x^2(x) dx$$
$$+ |k(\rho(b))|^2 + \left(1 + \frac{2(r+1)}{\rho_{\min}} F\right) |k(\rho(a))|^2 \quad (4.44)$$
$$+ \frac{\rho_{\max}}{\rho_{\min}} (r+4)^2 v^2(a) + (r+2) \frac{2\Gamma \rho_{\max}}{\rho_{\min}^2} \int_a^b \mu(\rho(x)) v_x^2(x) dx$$

Combining (4.44) with (4.30), we get:

$$V(a,b,\rho,v) \leq \left(c_2 + \frac{\tilde{K}(r+1)}{\rho_{\min}^2}\right) \int_a^b \rho^{-2}(x) P'(\rho(x)) \mu(\rho(x)) \rho_x^2(x) dx$$
$$+ \frac{R}{c_1}\left(2 + \frac{2(r+1)}{\rho_{\min}} F + \frac{\rho_{\max}}{\rho_{\min}}(r+4)^2\right) \min\left(((r+1)v(a) + k(\rho(a)))^2, v^2(a)\right)$$
$$+ \frac{1}{c_1}\left(2 + \frac{2(r+1)}{\rho_{\min}} F + \frac{\rho_{\max}}{\rho_{\min}}(r+4)^2\right) k(\rho(b))(P(\rho(b)) - P_{ext}) \quad (4.45)$$
$$+ \frac{1}{c_1}\left(2 + \frac{2(r+1)}{\rho_{\min}} F + \frac{\rho_{\max}}{\rho_{\min}}(r+4)^2\right) k(\rho(a))(P(\rho(a)) - P_{ext})$$
$$+ (r+2) \frac{2\Gamma \rho_{\max}}{\rho_{\min}^2} \int_a^b \mu(\rho(x)) v_x^2(x) dx$$

The existence of a constant $L > 0$ that depends only on $S \geq 0$, $r > 0$, $\rho^* > 0$ and satisfies (3.25) is a direct consequence of (4.45). The proof is complete. ◁

**Proof of Lemma 4:** Let $S \geq 0$ be given (arbitrary) and let $a,b \in \Re$ with $a < b$, $\rho \in C^0([a,b];(0,+\infty)) \cap H^1(a,b)$, $v \in C^0([a,b])$ with $\int_a^b \rho(x) dx = 1$ be given (arbitrary). Suppose that (3.11) holds. By virtue of Lemma 1, there exist constants $\rho_{\max} \geq \rho_{\min} > 0$ (depending only on $S \geq 0$, $r > 0$ and $\rho^* > 0$) such that (3.12) holds. Moreover, using the fact that $\int_a^b \rho(x) dx = 1$ in conjunction with (3.12), we obtain inequalities (4.33).

Using (3.7), (3.8), (3.9), (3.10) and the inequality $(x+y)^2 \leq 2x^2 + 2y^2$ that holds for all $x, y \in \Re$, we obtain the following inequality:

$$V(a,b,\rho,v) \leq \int_a^b \rho^{-3}(x)\mu^2(\rho(x))\rho_x^2(x)dx$$
$$+ |k(\rho(b))|^2 + |k(\rho(a))|^2 + (r+1)\int_a^b Q(\rho(x))dx \qquad (4.46)$$
$$+ (r+2)\frac{v^2(a)}{2} + (r+2)\frac{v^2(b)}{2} + \frac{r+2}{2}\int_a^b \rho(x)v^2(x)dx$$

Definitions (3.5), (3.6) imply that $Q''(\rho) = \rho^{-1}P'(\rho) > 0$, $k'(\rho) = \rho^{-1}\mu(\rho) > 0$ for all $\rho > 0$ and $Q(\rho^*) = Q'(\rho^*) = k(\rho^*) = 0$. Consequently, there exist constants $c_i > 0$, $i = 1,...,6$ (that depend only on $S \geq 0$, $r > 0$ and $\rho^* > 0$) such that the following inequalities hold:

$$c_1(\rho - \rho^*)^2 \leq Q(\rho) \leq c_2(\rho - \rho^*)^2, \text{ for all } \rho \in [\rho_{\min}, \rho_{\max}] \qquad (4.47)$$

$$c_3|\rho - \rho^*| \leq |k(\rho)| \leq c_4|\rho - \rho^*|, \text{ for all } \rho \in [\rho_{\min}, \rho_{\max}] \qquad (4.48)$$

$$c_5 \leq \rho^{-3}\mu^2(\rho) \leq c_6, \text{ for all } \rho \in [\rho_{\min}, \rho_{\max}] \qquad (4.49)$$

Using (4.47), (4.48), (4.49) in conjunction with (4.46) and (3.12), we obtain the following inequality:

$$V(a,b,\rho,v) \leq c_6 \|\rho_x\|_2^2 + c_2(r+1)\|\rho - \rho^*\chi_{[a,b]}\|_2^2 + \frac{r+2}{2}\rho_{\max}\|v\|_2^2$$
$$+ c_4^2|\rho(b) - \rho^*|^2 + c_4^2|\rho(a) - \rho^*|^2 + (r+2)\frac{v^2(a)}{2} + (r+2)\frac{v^2(b)}{2} \qquad (4.50)$$

Inequalities (4.33) imply the following estimate:

$$\|\rho - \rho^*\chi_{[a,b]}\|_2^2 \leq (b-a)\|\rho - \rho^*\chi_{[a,b]}\|_\infty^2 \leq \frac{1}{\rho_{\min}}\|\rho - \rho^*\chi_{[a,b]}\|_\infty^2 \qquad (4.51)$$

Combining (4.50) with (4.51) and using (4.33), we get:

$$V(a,b,\rho,v) \leq c_6\rho_{\max}(b-a)\|\rho_x\|_2^2 + \left(\frac{c_2(r+1)}{\rho_{\min}} + 2c_4^2\right)\|\rho - \rho^*\chi_{[a,b]}\|_\infty^2$$
$$+ \frac{\rho_{\max}(r+2)}{2\rho_{\min}}(b-a)^{-1}\|v\|_2^2 + (r+2)\frac{v^2(a)}{2} + (r+2)\frac{v^2(b)}{2} \qquad (4.52)$$

The first inequality (3.26) is a direct consequence of (4.52).

Using (3.7), (3.8), (3.9), (3.10) and the inequalities

$$v(x)\rho^{-1}(x)\mu(\rho(x))\rho_x(x) \geq -\frac{\varepsilon}{2}\rho(x)v^2(x) - \frac{1}{2\varepsilon}\rho^{-3}(x)\mu^2(\rho(x))\rho_x^2(x)$$

$$v(a)k(\rho(a)) \geq -\frac{\varepsilon}{2}v^2(a) - \frac{1}{2\varepsilon}|k(\rho(a))|^2$$

$$-v(b)k(\rho(b)) \geq -\frac{\varepsilon}{2}v^2(b) - \frac{1}{2\varepsilon}|k(\rho(b))|^2$$

which hold for all $\varepsilon > 0$, we obtain the following estimate for all $\varepsilon > 0$:

$$\begin{aligned}V(a,b,\rho,v) &\geq \frac{1}{2}\left(1-\varepsilon^{-1}\right)\int_a^b \rho^{-3}(x)\mu^2(\rho(x))\rho_x^2(x)dx \\ &+ \frac{1}{2}\left(1-\varepsilon^{-1}\right)|k(\rho(a))|^2 + \frac{1}{2}\left(1-\varepsilon^{-1}\right)|k(\rho(b))|^2 \\ &+ \frac{r+1-\varepsilon}{2}v^2(a) + \frac{r+1-\varepsilon}{2}v^2(b) + \frac{r+1-\varepsilon}{2}\int_a^b \rho(x)v^2(x)dx\end{aligned} \quad (4.53)$$

Selecting $\varepsilon = \dfrac{r+2}{2}$ and using (3.12), (4.49) we get:

$$V(a,b,\rho,v) \geq \frac{r}{2(r+2)}c_5\|\rho_x\|_2^2 + \frac{r}{4}v^2(a) + \frac{r}{4}v^2(b) + \frac{r}{4}\rho_{\min}\|v\|_2^2 + \frac{r}{2(r+2)}|k(\rho(a))|^2 \quad (4.54)$$

Exploiting the Cauchy-Schwarz inequality we get for all $x \in [a,b]$

$$\left|\rho(x) - \rho^*\right| = \left|\rho(a) - \rho^* + \int_a^x \rho_s(s)ds\right|$$

$$\leq |\rho(a) - \rho^*| + (x-a)^{1/2}\left(\int_a^x \rho_s^2(s)ds\right)^{1/2}$$

$$\leq |\rho(a) - \rho^*| + (b-a)^{1/2}\|\rho_x\|_2$$

from which we obtain the estimate (using the inequality $(x+y)^2 \leq 2x^2 + 2y^2$ that holds for all $x, y \in \Re$):

$$\|\rho_x\|_2^2 \geq \frac{1}{2(b-a)}\|\rho - \rho^*\chi_{[a,b]}\|_\infty^2 - \frac{1}{(b-a)}|\rho(a) - \rho^*|^2 \quad (4.55)$$

Using (4.48) and (4.54) we get:

$$\begin{aligned}V(a,b,\rho,v) &\geq \frac{r}{4(r+2)}\min\left(c_5, \frac{c_3^2}{\rho_{\max}}\right)\|\rho_x\|_2^2 + \frac{r}{4}v^2(a) + \frac{r}{4}v^2(b) \\ &+ \frac{r}{4}\rho_{\min}\|v\|_2^2 + \frac{r}{2(r+2)}c_3^2|\rho(a) - \rho^*|^2 + \frac{rc_5}{4(r+2)}\|\rho_x\|_2^2\end{aligned} \quad (4.56)$$

Combining (4.55) and (4.56), we obtain the following estimate:

$$V(a,b,\rho,v) \geq \frac{r}{8(r+2)(b-a)} \min\left(c_5, \frac{c_3^2}{\rho_{max}}\right) \|\rho - \rho^* \chi_{[a,b]}\|_\infty^2 + \frac{r}{4} v^2(a) + \frac{r}{4} v^2(b)$$
$$+ \frac{r}{4} \rho_{min} \|v\|_2^2 + \frac{rc_5}{4(r+2)} \|\rho_x\|_2^2 + \frac{r}{4(r+2)}\left(2c_3^2 - \frac{1}{(b-a)} \min\left(c_5, \frac{c_3^2}{\rho_{max}}\right)\right) |\rho(a) - \rho^*|^2 \quad (4.57)$$

Using inequalities (4.33) we obtain from (4.57) the following inequality

$$V(a,b,\rho,v) \geq \frac{r\rho_{min}}{8(r+2)} \min\left(c_5, \frac{c_3^2}{\rho_{max}}\right) \|\rho - \rho^* \chi_{[a,b]}\|_\infty^2 + \frac{r}{4} v^2(a) + \frac{r}{4} v^2(b)$$
$$+ \frac{r}{4} \rho_{min} \|v\|_2^2 + \frac{rc_5}{4(r+2)} \|\rho_x\|_2^2 + \frac{r}{4(r+2)}\left(2c_3^2 - \rho_{max} \min\left(c_5, \frac{c_3^2}{\rho_{max}}\right)\right) |\rho(a) - \rho^*|^2$$

which directly implies the estimate:

$$V(a,b,\rho,v) \geq \frac{r\rho_{min}}{8(r+2)} \min\left(c_5, \frac{c_3^2}{\rho_{max}}\right) \|\rho - \rho^* \chi_{[a,b]}\|_\infty^2$$
$$+ \frac{r}{4} v^2(a) + \frac{r}{4} v^2(b) + \frac{r}{4} \rho_{min} \|v\|_2^2 + \frac{rc_5}{4(r+2)} \|\rho_x\|_2^2 \quad (4.58)$$

Finally, using inequalities (4.33) in conjunction with (4.58) gives:

$$V(a,b,\rho,v) \geq \frac{r\rho_{min}}{8(r+2)} \min\left(c_5, \frac{c_3^2}{\rho_{max}}\right) \|\rho - \rho^* \chi_{[a,b]}\|_\infty^2$$
$$+ \frac{r}{4} v^2(a) + \frac{r}{4} v^2(b) + \frac{r\rho_{min}}{4\rho_{max}} (b-a)^{-1} \|v\|_2^2 + \frac{rc_5\rho_{min}}{4(r+2)} (b-a) \|\rho_x\|_2^2 \quad (4.59)$$

The second inequality (3.26) is a direct consequence of (4.52). The proof is complete. ◁

We end this section by providing the proofs of Theorem 1 and Theorem 2.

**Proof of Theorem 1:** Let $R, r > 0$, $S \geq 0$ be given. Lemma 3 and Lemma 4 imply that there exist constants $L > 0$, $G_2 \geq G_1 > 0$ that depend only on $S \geq 0$, $R, r > 0$ and $\rho^* > 0$ such that for every $a, b \in \Re$ with $a < b$, $\rho \in C^0([a,b];(0,+\infty)) \cap H^1(a,b)$, $v \in C^0([a,b]) \cap H^1(a,b)$ with $\int_a^b \rho(x)dx = 1$ and $V(a,b,\rho,v) \leq S$, inequalities (3.25), (3.26) hold.

Consider a classical solution for the PDE-ODE system (2.1)-(2.6), (3.18), i.e., consider functions $a, b \in C^1(\Re_+) \cap C^2((0,+\infty))$ satisfying $b(t) > a(t)$ for all $t \geq 0$, $\rho \in C^1(\bar{\Omega};(0,+\infty)) \cap C^2(\Omega)$, $v \in C^0(\bar{\Omega}) \cap C^1(\tilde{\Omega})$ with $v[t] \in C^2((a(t),b(t)))$ for each $t > 0$ and $\Omega = \bigcup_{t>0}(\{t\} \times (a(t),b(t)))$, $\tilde{\Omega} = \bigcup_{t>0}(\{t\} \times [a(t),b(t)])$ that satisfy equations (2.1)-(2.6) and (3.18). Assume that $V(a(0),b(0),\rho[0],v[0]) \leq S$.

Definition (3.7), equation (3.18) and Lemma 1 imply that the following equation holds for $t > 0$:

$$\frac{d}{dt}V\left(a(t),b(t),\rho[t],v[t]\right)=-\int_{a(t)}^{b(t)}\rho^{-2}(t,x)P'(\rho(t,x))\mu(\rho(t,x))\rho_x^2(t,x)dx$$
$$-k(\rho(t,b(t)))\left(P(\rho(t,b(t)))-P_{ext}\right)-k(\rho(t,a(t)))\left(P(\rho(t,a(t)))-P_{ext}\right) \quad (4.60)$$
$$-R\left((r+1)v(t,a(t))+k(\rho(t,a(t)))\right)^2 - r\int_{a(t)}^{b(t)}\mu(\rho(t,x))v_x^2(t,x)dx$$

Definition (3.5), equation (2.9) and the fact that $P'(\rho)>0$ for all $\rho>0$ imply that $k(\rho)\left(P(\rho)-P_{ext}\right)\geq 0$. Consequently, (4.60) implies the following inequality for all $t>0$:

$$\frac{d}{dt}V\left(a(t),b(t),\rho[t],v[t]\right)\leq 0 \quad (4.61)$$

Using (4.61) and the fact that the mapping $\Re_+ \ni t \to V\left(a(t),b(t),\rho[t],v[t]\right)$ is continuous, we conclude that

$$V\left(a(t),b(t),\rho[t],v[t]\right)\leq V\left(a(0),b(0),\rho[0],v[0]\right)\leq S, \text{ for all } t\geq 0 \quad (4.62)$$

Combining (3.25), (4.62) and (4.60) we get for all $t>0$:

$$\frac{d}{dt}V\left(a(t),b(t),\rho[t],v[t]\right)\leq -L^{-1}V\left(a(t),b(t),\rho[t],v[t]\right) \quad (4.63)$$

The differential inequality (4.63) and the fact that the mapping $\Re_+ \ni t \to V\left(a(t),b(t),\rho[t],v[t]\right)$ is continuous imply the following estimate:

$$V\left(a(t),b(t),\rho[t],v[t]\right)\leq \exp\left(-L^{-1}t\right)V\left(a(0),b(0),\rho[0],v[0]\right), \text{ for all } t\geq 0 \quad (4.64)$$

Define for all $t\geq 0$:

$$\Xi(t):=v^2(t,a(t))+v^2(t,b(t))+\left(\max_{a(t)\leq x\leq b(t)}\left(|\rho(t,x)-\rho^*|\right)\right)^2$$
$$+\frac{1}{b(t)-a(t)}\int_{a(t)}^{b(t)}v^2(t,x)dx+(b(t)-a(t))\int_{a(t)}^{b(t)}\rho_x^2(t,x)dx \quad (4.65)$$

Equation (3.20) implies the following inequalities for all $t\geq 0$:

$$\frac{1}{4}\Xi(t)\leq \left\|\Phi(a(t),b(t),\rho[t],v[t])-\left(\rho^*\chi_{[0,1]},0\chi_{[0,1]}\right)\right\|_X^2 \leq 4\Xi(t) \quad (4.66)$$

Combining (3.26), (4.62) and definition (4.65) we obtain for all $t\geq 0$:

$$G_1\Xi(t)\leq V(a(t),b(t),\rho[t],v[t])\leq G_2\Xi(t) \quad (4.67)$$

Combining (4.64), (4.66) and (4.67) we obtain estimate (3.19) with $M:=4\sqrt{\dfrac{G_2}{G_1}}$ and $\sigma:=\dfrac{1}{2L}$. The proof is complete. ◁

**Proof of Theorem 2:** Let $R > 0$, $S \geq 0$, $r \geq RK/2$ be given. Lemma 3 and Lemma 4 imply that there exist constants $L > 0$, $G_2 \geq G_1 > 0$ that depend only on $S \geq 0$, $R, r > 0$ and $\rho^* > 0$ such that for every $a, b \in \Re$ with $a < b$, $\rho \in C^0([a,b];(0,+\infty)) \cap H^1(a,b)$, $v \in C^0([a,b]) \cap H^1(a,b)$ with $\int_a^b \rho(x)dx = 1$ and $V(a,b,\rho,v) \leq S$, inequalities (3.25), (3.26) hold.

Consider a classical solution for the PDE-ODE system (2.1)-(2.6), (3.21), i.e., consider functions $a, b \in C^1(\Re_+) \cap C^2((0,+\infty))$ satisfying $b(t) > a(t)$ for all $t \geq 0$, $\rho \in C^1(\bar{\Omega};(0,+\infty)) \cap C^2(\Omega)$, $v \in C^0(\bar{\Omega}) \cap C^1(\tilde{\Omega})$ with $v[t] \in C^2((a(t),b(t)))$ for each $t > 0$ and $\Omega = \bigcup_{t>0}(\{t\} \times (a(t),b(t)))$, $\tilde{\Omega} = \bigcup_{t>0}(\{t\} \times [a(t),b(t)])$ that satisfy equations (2.1)-(2.6) and (3.21). Assume that $V(a(0),b(0),\rho[0],v[0]) \leq S$.

Definition (3.7), equation (3.21) and Lemma 1 imply that the following equation holds for $t > 0$:

$$\frac{d}{dt}V(a(t),b(t),\rho[t],v[t]) = -\int_{a(t)}^{b(t)} \rho^{-2}(t,x)P'(\rho(t,x))\mu(\rho(t,x))\rho_x^2(t,x)dx$$

$$-k(\rho(t,b(t)))\left(P(\rho(t,b(t))) - P_{ext}\right) - k(\rho(t,a(t)))\left(P(\rho(t,a(t))) - P_{ext}\right) \quad (4.68)$$

$$-R(r+1)v^2(t,a(t)) - Rk(\rho(t,a(t)))v(t,a(t)) - r\int_{a(t)}^{b(t)} \mu(\rho(t,x))v_x^2(t,x)dx$$

Definition (3.5), equation (2.9) and the fact that $P'(\rho) > 0$ for all $\rho > 0$ imply that $k(\rho)(P(\rho) - P_{ext}) \geq 0$. Consequently,

$$k(\rho(t,a(t)))\left(P(\rho(t,a(t))) - P_{ext}\right) = |k(\rho(t,a(t)))||P(\rho(t,a(t))) - P_{ext}| \quad (4.69)$$

Using the inequality $-k(\rho(t,a(t)))v(t,a(t)) \leq \frac{1}{2RK}|k(\rho(t,a(t)))|^2 + \frac{RK}{2}v^2(t,a(t))$, we obtain from (4.68) and (4.69) for $t > 0$:

$$\frac{d}{dt}V(a(t),b(t),\rho[t],v[t]) \leq -\int_{a(t)}^{b(t)} \rho^{-2}(t,x)P'(\rho(t,x))\mu(\rho(t,x))\rho_x^2(t,x)dx$$

$$-k(\rho(t,b(t)))\left(P(\rho(t,b(t))) - P_{ext}\right) - R\left(r+1-\frac{RK}{2}\right)v^2(t,a(t)) \quad (4.70)$$

$$-|k(\rho(t,a(t)))|\left(|P(\rho(t,a(t))) - P_{ext}| - \frac{1}{2K}|k(\rho(t,a(t)))|\right) - r\int_{a(t)}^{b(t)} \mu(\rho(t,x))v_x^2(t,x)dx$$

Using (3.22), (2.9), (4.69) and the fact that $r \geq RK/2$, we obtain from (4.70) for $t > 0$:

$$\frac{d}{dt}V(a(t),b(t),\rho[t],v[t]) \leq -\int_{a(t)}^{b(t)} \rho^{-2}(t,x)P'(\rho(t,x))\mu(\rho(t,x))\rho_x^2(t,x)dx$$

$$-k(\rho(t,b(t)))\left(P(\rho(t,b(t))) - P_{ext}\right) - \frac{1}{2}k(\rho(t,a(t)))\left(P(\rho(t,a(t))) - P_{ext}\right) \quad (4.71)$$

$$-Rv^2(t,a(t)) - r\int_{a(t)}^{b(t)} \mu(\rho(t,x))v_x^2(t,x)dx$$

From this point on, the proof is completely the same with the proof of Theorem 1. The proof is complete. ◁

# 5. Concluding Remarks

The present paper provides globally stabilizing feedback laws for a two-piston problem that is based on the nonlinear Navier-Stokes equations for a compressible fluid. We showed that the CLF methodology can be applied successfully by constructing an appropriate CLF for the system. The proposed feedback laws have minimal measurement requirements and it was shown that all assumptions are valid for a wide class of gases (which includes the ideal gas).

However, the main results leave some open problems:

1) The results were applied to classical solutions. It is of interest to relax this to weak solutions. It is an open problem to show that the proposed feedback laws preserve their strong stabilizing role in the case of weak solutions.

2) The assumption of isentropic (or barotropic) gas flow may also be relaxed. It is an open problem to show that the proposed feedback laws can achieve global stabilization when the energy balance PDE is added to the model.

It is not only interesting to examine but, in fact, essential to understand the relation between the present work and the results in [19] and [21] for the viscous Burgers' equation. The viscous Burgers PDE has been used for decades as a simplified surrogate in the study of incompressible fluid flows, especially in one dimension.

Paper [21], which, incidentally, is the first paper on applying ODE backstepping in PDE control, following the work [5] by Coron and d'Andrea-Novel outside of a fluid context, considered inputs acting through integrators, in (1.2) and (1.3) in [21], just as in the case of an actuated piston in (2.5) of the present paper, with the second big difference (besides the constant density) being that in the present paper the boundary's location is not fixed. The control laws (2.8), (2.9) in [21] are far more complex than the feedback (3.18) here (and employ even the boundary values of the second spatial derivative of velocity). The greater simplicity of the control law (3.18) stems, ironically, from exploiting the greater complexity of the model (2.1)-(2.6), including all of its extra features (non-constant pressure, moving boundaries, and conservation of mass).

A comparison with the flux-actuated viscous (and incompressible) Burgers' equation in (2.1), (2.4), (2.5) of [19], which is a simpler problem than the problem in [21], and in a way more distant than the problem in this paper, is also interesting. Starting from the end result of the feedback design in [19], given by the quadratic feedback in (3.3), (3.4), or the cubic feedback (3.7), (3.8), one cannot help but notice the uncanny similarity with the feedback (3.18) here, and even more with the linear velocity feedback (3.21). This should not be surprising. Note that, if the left piston were heavy, a singular perturbation argument would result in the left side of (2.5) being zero, and in the actuation of the system being of the Neumann type in velocity, just as in (2.4), (2.5) of [19]. The reason why the quadratic or cubic terms are needed in (3.3), (3.4) or (3.7), (3.8), respectively, in [19], while such terms are not seen in (3.21), is that, again, we exploit the additional complexities of the compressible moving-boundary problem here, whereas, with those additional physical complexities missing in the viscous Burgers' model in [19], nonlinear feedback of velocity is needed for Neumann-actuated stabilization.

It is important to recognize that none of these physical systems - neither the one here nor the ones in [19,21] - are open-loop unstable. However, all these systems lack asymptotic stability, which is most easily seen in [19] where the linearization has an eigenvalue at the origin. Hence, all the three designs, including the one here, endeavor to endow a neutrally stable highly nonlinear system with asymptotic stability, without destroying the global character of the neutral stability of the open-loop system.

Finally, since the present paper deals with two coupled PDEs, with ODEs on each end of their 1-D domain, and with control applied to only one of these two PDEs, it should not be overlooked that such a "sandwich" structure originates in [35], however, without moving boundaries, nonlinearities, and viscosity, but with potential open-loop instability in both the PDEs and the ODEs.

**References**

[1] Annaswamy, A. M. and A. F. Ghoniem, "Active Control in Combustion Systems," *Control Systems Magazine*, 15, 1995, 49–63.

[2] Auriol, J., G. A. de Andrade, R. Vazquez, "A Differential-Delay Estimator for Thermoacoustic Oscillations in a Rijke Tube Using In-Domain Pressure Measurements", *Proceedings of the 2020 59th IEEE Conference on Decision and Control (CDC)*, 2020, 4417-4422.

[3] Banaszuk, A., K. B. Ariyur, and M. Krstic, "An Adaptive Algorithm for Control of Combustion Instability", *Automatica*, 40, 2004, 1965–1972.

[4] Buisson-Fenet, M., S. Koga and M. Krstic, "Control of Piston Position in Inviscid Gas by Bilateral Boundary Actuation", *Proceedings of the 2018 57th IEEE Conference on Decision and Control (CDC)*, 2018, 5622-5627.

[5] Coron, J.-M. and B. d'Andrea-Novel, "Stabilization of a Rotating Body Beam without Damping", *IEEE Transactions on Automatic Control*, 43, 1998, 608–618.

[6] Coron, J.-M., *Control and Nonlinearity*, Mathematical Surveys and Monographs, Volume 136, American Mathematical Society, 2007.

[7] de Andrade, G. A., R. Vazquez, and D. J. Pagano, "Backstepping Stabilization of a Linearized ODE-PDE Rijke Tube Model", *Automatica*, 96, 2018, 98-109.

[8] de Andrade, G.A., R. Vazquez and D. J. Pagano, "Backstepping-Based Estimation of Thermoacoustic Oscillations in a Rijke Tube With Experimental Validation", *IEEE Transactions on Automatic Control*, 65, 2020, 5336-5343.

[9] Dick, M., M. Gugat and G. Leugering, "A Strict $H^1$-Lyapunov Function and Feedback Stabilization for the Isothermal Euler Equations With Friction", *Numerical Algebra, Control and Optimization*, 1, 2011, 225-244.

[10] Epperlein, J. P., B. Bamieh and J. Astrom, "Thermoacoustics and the Rijke Tube: Experiments, Identification and Modeling", *Control Systems Magazine*, 35, 2015, 57-77.

[11] Feireisl, E., *Dynamics of Viscous Compressible Fluids*, Oxford University Press, Oxford, 2004.

[12] Feireisl, E., V. Macha, S. Necasova and M. Tucsnak, "Analysis of the Adiabatic Piston Problem via Method of Continuum Mechanics", *Annales de l'Institut Henri Poincaré* C, *Analyse non Linéaire*, 35, 2018, 1377-1408.

[13] Guenther, R. B., and J. W. Lee, *Partial Differential Equations of Mathematical Physics and Integral Equations*, Dover, 1996.

[14] Gugat, M. and M. Herty, "Existence of Classical Solutions and Feedback Stabilization for the Flow in Gas Networks", *ESAIM: Control, Optimisation and Calculus of Variations*, 17, 2011, 28-51.

[15] Karafyllis, I. and M. Krstic, "Global Stabilization of a Class of Nonlinear Reaction-Diffusion PDEs by Boundary Feedback", *SIAM Journal on Control and Optimization*, 57, 2019, 3723-3748.

[16] Karafyllis, I., "Lyapunov-Based Boundary Feedback Design For Parabolic PDEs", *International Journal of Control*, 94(5), 2021, 1247–1260.

[17] Kazhikhov, A. V. and V. V. Shelukhin, "Unique Global Solution with Respect to Time of Initial-Boundary Value Problems for One-Dimensional Equations of a Viscous Gas", *Journal of Applied Mathematics and Mechanics*, 41, 1977, 273-282, 1977.

[18] Koga, S. and M. Krstic, *Materials Phase Change PDE Control & Estimation*, Birkhäuser, 2020.

[19] Krstic, M., "On Global Stabilization of Burgers' Equation by Boundary Control", *Systems and Control Letters*, 37, 1999, 123-142.

[20] Lions, P.-L., *Mathematical Topics in Fluid Dynamics, Vol.2, Compressible Models*, Oxford Science Publication, Oxford, 1998.

[21] Liu, W.-J. and M. Krstic, "Backstepping Boundary Control of Burgers' Equation with Actuator Dynamics," *Systems & Control Letters*, 41, 2000, 291-303.

[22] Liu, Y., T. Takahashi and M. Tucsnak, "Single Input Controllability of a Simplified Fluid-Structure Interaction Model", *ESAIM: Control, Optimization and Calculus of Variations*, 19, 2013, 10–42.

[23] Maity, D., T. Takahashi and M. Tucsnak, "Analysis of a System Modelling the Motion of a Piston in a Viscous Gas", *Journal of Mathematical Fluid Mechanics*, 19, 2017, 551–579.

[24] Moore, W. J., *Physical Chemistry*, 5th Edition, Longman, 1998.

[25] Nash, J., "Le Problème de Cauchy pour les Equations Différentielles d'un Fluide Général", *Bulletin de la Société Mathématique de France*, 90, 1962, 487-497.

[26] Plotnikov, P. I. and J. Sokolowski, "Boundary control of the Motion of a Heavy Piston in Viscous Gas", *SIAM Journal on Control and Optimization*, 57, 2019, 3166-3192.

[27] Shelukhin, V. V., "The Unique Solvability of the Problem of Motion of a Piston in a Viscous Gas", *Dinamika Sploshnoi Sredy*, 31, 1977, 132–150.

[28] Shelukhin V. V., "Stabilization of the Solution of a Model Problem on the Motion of a Piston in a Viscous Gas", *Dinamika Sploshnoi Sredy*, 173, 1978, 134–146.

[29] Shelukhin, V. V., "Motion with a Contact Discontinuity in a Viscous Heat Conducting Gas", *Dinamika Sploshnoi Sredy*, 57, 1982, 131–152.

[30] Smoller, J., *Shock Waves and Reaction-Diffusion Equations*, 2nd Edition, Springer-Verlag, New York, 1994.

[31] Smyshlyaev, A. and M. Krstic, "Closed-Form Boundary State Feedbacks for a Class of 1-D Partial Integro-Differential Equations", *IEEE Transactions on Automatic Control*, 49, 2004, 2185-2202.

[32] Smyshlyaev, A. and M. Krstic, *Adaptive Control of Parabolic PDEs*, Princeton University Press, 2010.

[33] Temam, R. M. and A. M. Miranville, *Mathematical Modeling in Continuum Mechanics*, 2nd Edition, Cambridge University Press, 2005.

[34] Vazquez, R. and M. Krstic, *Control of Turbulent and Magnetohydrodynamic Channel Flows*, Birkhäuser, 2007.

[35] Wang, J., M. Krstic and Y. Pi, "Control of a 2×2 Coupled Linear Hyperbolic System Sandwiched Between 2 ODEs", *International Journal of Robust and Nonlinear Control*, 2018, 1–30.